\documentclass[10pt]{article}
\usepackage{latexsym}
\usepackage{amsmath, amsthm, amsfonts}
\usepackage{graphics}
\usepackage{graphicx}
\usepackage{color}
\usepackage[table,xcdraw]{xcolor}
\usepackage{float}
\usepackage{array}
\usepackage{hyperref}

\newcommand{\abs}[1]{\left\lvert#1\right\rvert}
\newtheorem{theorem}{Theorem}[section]

\def\smallskip{\addvspace{\smallskipamount}}
\def\medskip{\addvspace{\medskipamount}}
\def\bigskip{\addvspace{\bigskipamount}}
\def\makefootline{\baselineskip=24pt \line{\the\footline}}
\def\pagecontents{\ifvoid\topins\else\unvbox\topins\fi
   \dimen@=\dp255 \unvbox255
   \ifvoid\footins\else
      \vskip\skip\footins \footnoterule \unvbox\footins\fi
     \ifr@ggedbottom \kern-\dimen@ \vfil \fi}
\def\footnoterule{\kern-3pt\hrule width 2truein \kern 2.6pt}

\begin{document}

\title{\textbf{COMPUTATIONAL COMPLEX DYNAMICS OF THE DISCRETE LORENZ SYSTEM}}

\author{Sk. Sarif Hassan\\
  \small {\emph{Department of Mathematics}}\\
\small {\emph{College of Engineering Studies}} \\
  \small {\emph{University of Petroleum and Energy Studies}}\\
  \small {\emph{Bidholi, Dehradun, India}}\\
  \small Email: {\texttt{\textcolor[rgb]{0.00,0.07,1.00}{s.hassan@ddn.upes.ac.in}}}\\
}

\maketitle
\begin{abstract}
\noindent The dynamics of the classical Lorenz system is well studied in $1963$ by E. N. Lorenz. Later on, there have been an extensive studies on the classical Lorenz system with the complex variables and the discrete time Lorenz system with real variables. To the best of knowledge of the author, so far there is no study on discrete time Lorenz system in complex variables. In this article, an attempt has been made to observe and understand the discrete dynamics of the Lorenz system with complex variables. This study compares the discrete dynamics of the Lorenz system with complex variables to that of the classical Lorenz system involving real and complex variables.    
\end{abstract}

\begin{flushleft}\footnotesize
{\textbf{Keywords:} Discrete Complex Lorenz System, Hyper Periodic, Transient chaos, Periodic \& Chaos. \\

}.
\end{flushleft}

\section{Introduction}
 The concept of chaos was first introduced the Lorenz system by \emph{E. Lorenz} \cite{EL1} and \cite{EL2} and it is the prototypical example of sensitive dependence on initial conditions (the butterfly effect) \cite{SI}. As time progresses, interest in the complex dynamic behaviors of nonlinear systems increased, due to their potential applications in different fields, such as detecting changes of biological signals (mostly EEG) in different abnormalities, data and image encryption, studying sunspot cycles, lasers and so on. Consequently, a large number of novel systems have been developed based on the original Lorenz system viz. R\"{o}ssler system, Chen system and so on \cite{C-J-C}, \cite{GX}, \cite{OR}, \cite{JCS1}, \cite{GCT}, \cite{JG}, \cite{JGS}, \cite{JGD}, \cite{JCS2} and \cite{JCS3}. It is interesting to note that in $2000$, there are eighteen challenging mathematical problems for the twenty-first century as introduced by a mathematician \emph{Smale} \cite{Smale}, in which the Lorenz system attractor was the fourteenth problem. In this regard, one of the main concern was about existence of the strange attractor which also has been proved by \emph{Stewart} \cite{SI}. \\

\noindent
The Lorenz system is described as \cite{EL1}:

\begin{equation}
\frac{\partial x}{\partial t}=a(y-x),
\frac{\partial y}{\partial t}=rx-y-xz,
\frac{\partial z}{\partial t}=xy-bz
\label{equation:total-equationA}
\end{equation}

\noindent
Some key features of this system are as follows \cite{WJ}:

\begin{itemize}
  \item It is an autonomous system, which means that time does not explicitly appear on the right hand side of the Eq.(1).
  \item The equations involve only first order time derivatives, so the evolution depends only on the values of $x$, $y$, and $z$ at the time.
  \item Due to the terms of $xz$ and $xy$ in the second and third equations in Eq.(1), the system is \emph{non-linear}.
  \item The system is \emph{dissipative} when the following inequality holds: $\frac{\partial^2 x}{\partial t^2}+\frac{\partial^2 y}{\partial t^2}+\frac{\partial^2 z}{\partial t^2}=-(a+b+1)<0$\\
      Since parameters $a$ and $b$, denoting the physical characteristics of the air flow, are positive, the inequality always holds and, thus, solutions are bounded.

  \item The system is \emph{symmetric}, with respect to the $z$ axis, which means it is invariant for the coordinate transformation: $(x, y, z) \rightarrow (-x,-y,z)$.
\end{itemize}

\noindent
The main observation by \emph{Lorenz} was that the Lorenz system exhibits the sensitive dependence on initial conditions-very small changes in initial conditions can make very large differences in long-term behavior. The trajectory $(x(t), y(t), z(t))$, when $r = 28$, approaches a \emph{strange attractor}, that is, an attractor that is not an equilibrium, nor a cycle, nor a finite graph. The same Lorenz system Eq.(\ref{equation:total-equationA}) has been seen as a discrete time system and studied the dynamics in \cite{WJ} with a consideration that the parameters and variables are real numbers. \\

\noindent
In this article, we consider the discrete time approximated Lorenz system in the complex plane and we would like to compare the dynamics which has been already seen by others in the real and complex variables of the Lorenz system \cite{MH}, \cite{GSM}, \cite{JWC}, \cite{AHB}, \cite{GMT}, \cite{GM}, \cite{GP}, \cite{JCS4}, \cite{WWM} and \cite{AA}. \\\\
\noindent
The discrete time Lorenz system is the following:

\begin{equation}
\frac{x_{k+1}-x_k}{dt}=a(y_k-x_k),
\frac{y_{k+1}-y_k}{dt}=-x_kz_k+rx_k-y_k,
\frac{z_{k+1}-z_k}{dt}=x_ky_k-bz_k
\label{equation:total-equationB}
\end{equation}

\noindent
Therefore the iterative system of Lorenz system is from the Eq.(\ref{equation:total-equationB}):
\begin{eqnarray}
  x_{k+1} &=& x_k+a(y_k-x_k)dt \\
  y_{k+1} &=& y_k+(-x_k z_k+rx_k-y_k)dt \\
  z_{k+1} &=& z_k+(x_ky_k-b z_k)dt
  \label{eqnarray:total-equationC}
\end{eqnarray}
\noindent
Here all the parameters $a$ (Prandtl number), $b$, $r$ (Rayleigh
number) and the initial values $x_0, y_0$ and $z_0$ are complex variables. Further in due course, we designate the Lorenz system Eq.$(3,4,5)$ as discrete complex Lorenz system. \\

\noindent
The key dynamics what we have achieved in the complex plane are summarized here \dots \\

\noindent
There are complex parameters $r$, $a$ and $b$ with complex initial values for which the discrete complex Lorenz system Eq.$(3,4,5)$ has the following kind of solutions/trajectories as found in this present computational study.
\begin{itemize}
  \item There are three fixed points of the Lorenz systems Eq.$(3,4,5)$ and there are certain parameters $r, a,$ and $b$ (examples are given) such that the fixed points are stable (sink).
  \item The system Eq.$(3,4,5)$ possesses higher order periodic solutions too (few examples are given).
  \item The system Eq.$(3,4,5)$ has chaotic and transient chaotic solutions.
  \item The system Eq.$(3,4,5)$ has single and double coexisting chaotic attractors and existence of the chaotic attractors has been assured through examples.
  \item A comparison have been made with other existing real and complex classical Lorenz systems.
\end{itemize}

\noindent
In the following sections, a detail dynamics of the discrete complex Lorenz system Eq.$(3,4,5)$ is characterized and compared with the existing real and complex classical Lorenz systems \cite{QMA}, \cite{JE}, \cite{TN}, \cite{TL} and \cite{HR}.

\section{Stability of the Fixed Points}

Without loss of generality, in the rest of the article, all the initial values are taken from unit disk $B(0,1)\subset\mathbb{C}$.
\noindent
The fixed points of the Lorenz system Eq.$(3,4,5)$ are the solutions of the system of equations:
\[
\bar{x}=\bar{x}+a(\bar{y}-\bar{x})dt; \bar{y}=\bar{y}+(-\bar{x}\bar{z}+r\bar{x}-\bar{y})dt; \bar{z}=\bar{z}+(\bar{x}\bar{y}-b\bar{z})dt
\]
\noindent
Consequently, the system Eq.$(3,4,5)$ has the unique three fixed points \dots \\ \\ $(\bar{x},\bar{y},\bar{z})_{1}=(0,0,0),$ $(\bar{x},\bar{y},\bar{z})_{2}=(-\sqrt{b} \sqrt{-1+r}, -\sqrt{b} \sqrt{-1+r}, -1+r)$ and $(\bar{x},\bar{y},\bar{z})_{3}=(\sqrt{b} \sqrt{-1+r}, \sqrt{b} \sqrt{-1+r}, -1+r)$.\\ \\
\noindent
It is nice to note that in classical real Lorenz system, if $r<1$ then there is only one equilibrium point, which is at the origin $(0,0,0)$. All orbits converge to $(0,0,0)$, which is a global attractor, when $r<1$. A pitchfork bifurcation occurs at $r=1$, and when $r>1$ then other two fixed points are stable only if $r< a\frac{a+b+3}{a-b-1}$ which can hold only for positive $r$ and if $a>b+1$.   \\

\noindent
The stability of the fixed points depends on the eigenvalues of the jacobian of the system of equations.

\begin{theorem}
The fixed points $(\bar{x},\bar{y},\bar{z})_{i}$ for $i=1,2$ and $3$ of the Eq.$(3,4,5)$ is \dots \\
\emph{stable/sink/attractor} if the real parts of the eigenvalues of the jacobian at the fixed point are positive.\\
\emph{unstable} if the real parts of the eigenvalues of the jacobian at the fixed point are positive. \\
\emph{non-hyperbolic} if the real parts of the eigenvalues of the jacobian at the fixed point are zero.
\end{theorem}
\noindent
The jacobian of the system is $J=\left(
                                   \begin{array}{ccc}
                                     \frac{\partial f}{\partial x} & \frac{\partial f}{\partial y} & \frac{\partial f}{\partial z} \\
                                     \frac{\partial g}{\partial x} & \frac{\partial g}{\partial y} & \frac{\partial g}{\partial z} \\
                                     \frac{\partial h}{\partial x} & \frac{\partial h}{\partial y} & \frac{\partial h}{\partial z} \\
                                   \end{array}
                                 \right)
$ where $f(x,y,z)=x+a(y-x)dt$; $g(x,y,z)=y+(-xz+rx-y)dt$ and $h(x,y,z)=z+(xy-bz)dt$.\\
\noindent
So the jacobian becomes $J(x,y,z)=\left(
                                    \begin{array}{ccc}
                                      1+a dt & -a dt & 0 \\
                                      (r+z)dt  & 1-dt & x dt  \\
                                      y dt & x dt & 1-b dt \\
                                    \end{array}
                                  \right)
$

\noindent
Here the jacobian $J$ at the fixed point $(\bar{x},\bar{y},\bar{z})_{1}=(0,0,0)$ is $J(0,0,0)=\left(
                \begin{array}{ccc}
                  1+adt & -adt & 0 \\
                  rdt & 1-dt & 0 \\
                  0 & 0 & 1-bdt \\
                \end{array}
              \right)$ \\ \\

\noindent
The eigenvalues of $J(0,0,0)$ are \dots \\ $\left\{1-b \text{dt},\frac{1}{2} \left(2-\text{dt}+a \text{dt}-\text{dt} \sqrt{1+2 a+a^2-4 a r}\right),\frac{1}{2} \left(2-\text{dt}+a \text{dt}+\text{dt} \sqrt{1+2 a+a^2-4 a r}\right)\right\}$. The fixed point $(\bar{x},\bar{y},\bar{z})_{1}=(0,0,0)$ is stable if the real parts of the eigenvalues of the Jacobian at the $(\bar{x},\bar{y},\bar{z})_{1}$ are positive. \\\\
\noindent
Similarly, the other eigenvalues can be archived at the fixed points at the jacobian \emph{J}. Here we omit the detail as the general form of the eigenvalues are complicated in looking. Let us have a few examples of stability of these fixed points.\\

\noindent
Consider the parameters $a=10$, $b=\frac{8}{3}$ and $r=5$, and $dt=0.005$ in the Lorenz system Eq.$(3,4,5)$. In this case, the fixed points are $(\bar{x},\bar{y},\bar{z})_{1}=(0,0,0)$, $(\bar{x},\bar{y},\bar{z})_{2}=(-4,-4,6)$ and $(\bar{x},\bar{y},\bar{z})_{3}=(4,4,6)$. For different initial values (a set of five initial values), we generate the trajectories as shown in Fig. $1$.

\begin{figure}[H]
      \centering

      \resizebox{10cm}{!}
      {
      \begin{tabular}{c c}
      \includegraphics [scale=6]{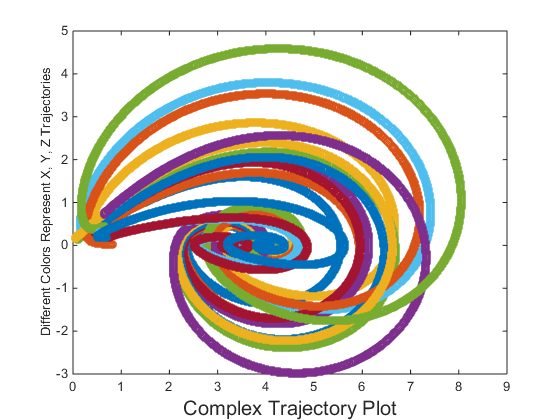}
      \includegraphics [scale=6]{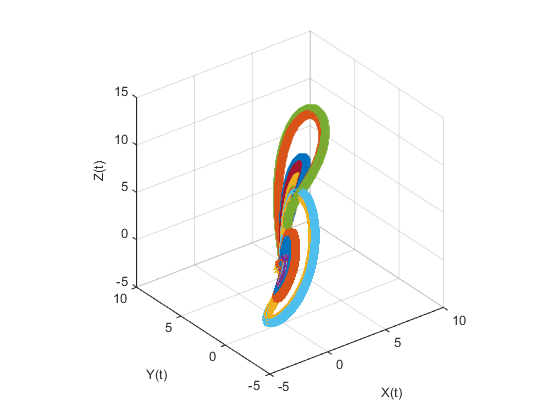}
      \end{tabular}
      }
\caption{Trajectory Plots (Left: Five X, Y, Z Complex trajectories, Right: Five trajectories in 3 dimension.)}
      \begin{center}

      \end{center}
      \end{figure}

\noindent
It is found that the eigenvalues of the jacobian \emph{J} at the fixed points $(\bar{x},\bar{y},\bar{z})_{2}=(-4,-4,6)$ and $(\bar{x},\bar{y},\bar{z})_{3}=(4,4,6)$ are  $0.977399,1.02713\, -0.0402936 i$ and $1.02713\, +0.0402936 i$. Hence the real part of all the eigenvalues are positive and therefore the fixed points $(\bar{x},\bar{y},\bar{z})_{2}=(-4,-4,6)$ and $(\bar{x},\bar{y},\bar{z})_{3}=(4,4,6)$ are stable (sink) which is evident from the Fig.$1$ too.\\
\noindent
It is noted that $\abs{a}>\abs{b+1}$ and $\abs{r}>1$ which is naturally the condition for local stability at the fixed points as it is observed in the case of classical real Lorenz system.  \\ \\
\noindent
Consider the parameters $a=10$, $b=\frac{8}{3}$ and $r=4+9i$, and $dt=0.0005$ in the Lorenz system $(3,4,5)$. In this case, the fixed points are $(\bar{x},\bar{y},\bar{z})_{1}=(0,0,0)$, $(\bar{x},\bar{y},\bar{z})_{2}=(-4.51599 - 2.65722i,-4.51599 - 2.65722i,5+9i)$ and $(\bar{x},\bar{y},\bar{z})_{3}=(4.51599 + 2.65722i,4.51599 + 2.65722i,5+9i)$. For a set of five initial values, we generate the trajectories as shown in Fig.$2$.

\begin{figure}[H]
      \centering

      \resizebox{11cm}{!}
      {
      \begin{tabular}{c c}
      \includegraphics [scale=8]{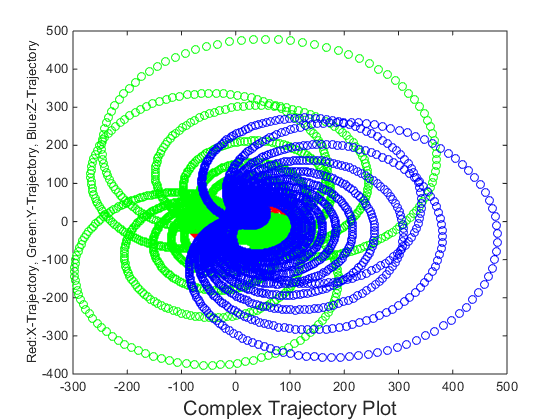}
      \includegraphics [scale=8]{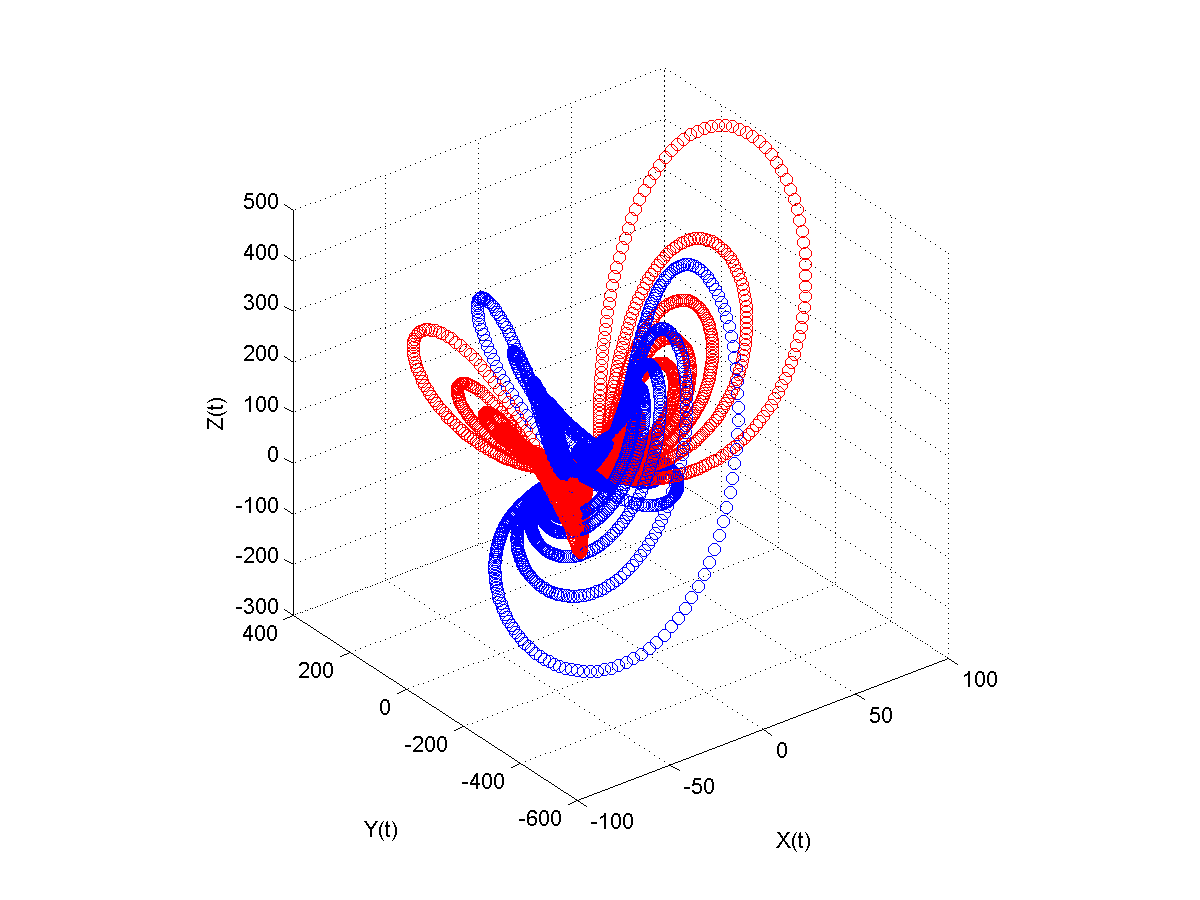}\\
      \includegraphics [scale=8]{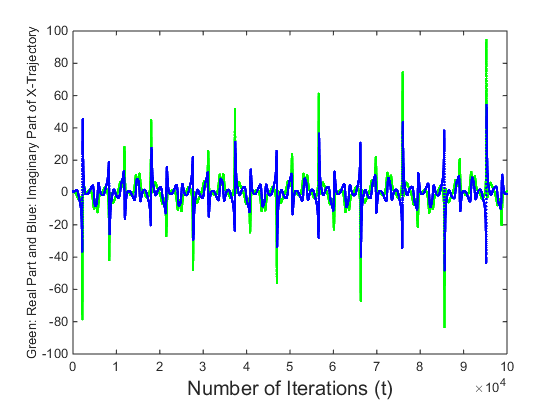}
      \includegraphics [scale=8]{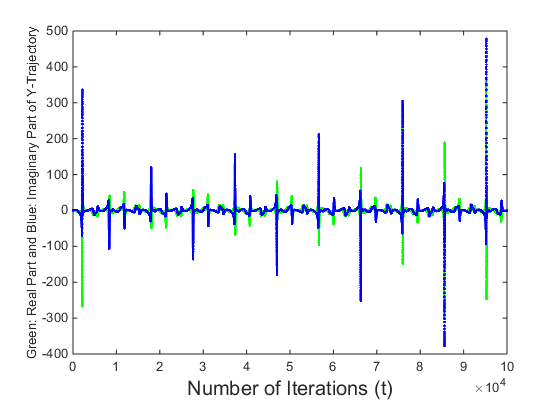}
      \includegraphics [scale=8]{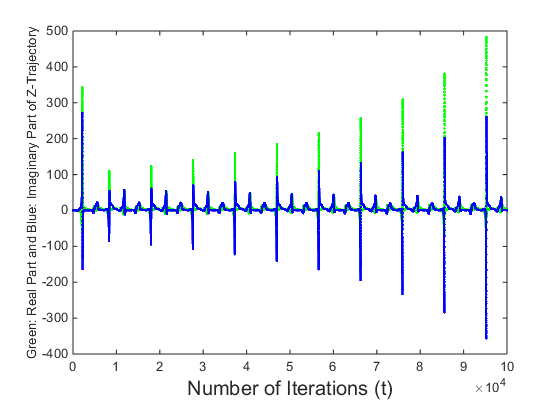}\\

            \end{tabular}
      }
\caption{Trajectory Plots (Top: A complex trajectory (Left) and Trajectory in 3 dimension (Right) and in Bottom: Real and Imaginary Part of X, Y, Z trajectories)}
      \begin{center}

      \end{center}
      \end{figure}

\noindent
It is found that the eigenvalues of the jacobian \emph{J} at the fixed points $(\bar{x},\bar{y},\bar{z})_{2}=(-4,-4,6)$ and $(\bar{x},\bar{y},\bar{z})_{3}=(4,4,6)$ are  $0.977399,1.02713\, -0.0402936 i$ and $1.02713\, +0.0402936 i$ are  ${0.973849 - 0.00410134i, 0.990298 + 0.0554783i, 1.06752 - 0.0513769i}$. Hence the real part of all the eigenvalues are positive and therefore the fixed points $(\bar{x},\bar{y},\bar{z})_{2}=(-4,-4,6)$ and $(\bar{x},\bar{y},\bar{z})_{3}=(4,4,6)$ are stable (sink) which is evident from the Fig.$2$ too. For this set of parameters, it is noted that $\abs{r}>1$ and $\abs{r}<\abs{a\frac{a+b+3}{a-b-1}}$ is holding well in discrete complex Lorenz system Eq.$(3,4,5)$.

\section{Dynamics in Discrete Complex Lorenz System}
Here we shall take different complex parameters and initial vales and run the dynamical system to observe the trajectory behaviors. Further we shall compare the present results with same of the classical discrete Lorenz systems. Here we first consider arbitrary complex initial values and the parameters are set as $a=10$ and $b=\frac{8}{3}$ and $dt=0.0005$.

\subsection{Comparison of Dynamics of the Discrete Real Lorenz System and Discrete Complex Lorenz System}
Considering the parameters $a=10$, $b=\frac{8}{3}$ and $dt=0.0005$ with different choices of the control parameter $r$, what happens in dynamics of the Lorenz system Eq.$(3,4,5)$ is our aim to understand.

\begin{table}[H]

\begin{tabular}{| m{1cm}||m{3cm}| |m{5cm}| |m{5cm}|}
\hline \centering \textbf{Serial No} &
\begin{center}
\textbf{Control Parameter $r$}
\end{center}
 &
\begin{center}
\textbf{Behavior in Discrete Complex Lorenz System}
\end{center}
&
\begin{center}
\textbf{Behavior in Discrete Real Lorenz System}
\end{center}
\\
\hline \centering 1 &
\centering $10$
&
\begin{center}
Converges to one of the sinks $(4.899,4.8990,9)$
\end{center}
&
\begin{center}
Converges to one of the sinks $(4.899,4.8990,9)$

\end{center}\\
\hline
\centering 2 &
\centering $18$
&
\centering  Converges to one of the sinks $(-6.74,-6.74,17)$

&
\begin{center}
Converges to one of the sinks $(-6.74,-6.74,17)$ (Transient Chaos)

\end{center}
\\
\hline
\centering 3 &
\centering $22.35$
&
\begin{center}
 Converges to one of the sinks $(-7.55 ,-7.55,21.35)$
\end{center}
&
\begin{center}
Chaotic
\end{center}
\\
\hline
\centering 4 &
\centering $26$
&
\begin{center}
Chaotic (Largest exponent is positive)
\end{center}
&
\begin{center}
Chaotic
\end{center}
\\
\hline
\centering 5 &
\centering $28$
&
\begin{center}
Transient Chaotic  (One positive Lyapunav exponent.)
\end{center}
&
\begin{center}
Chaotic
\end{center}
\\
\hline

\end{tabular}
\caption{Dynamics of complex Lorenz System for different values of $r$.}
\label{Table:1}
\end{table}
\noindent
In the discrete real Lorenz system, it is found that the system possesses chaotic behavior when $r \geq 22.35$ with $a=10$ and $b=\frac{8}{3}$. Also for $r=18$, the trajectories in the discrete Lorenz system will appear to be chaotic for an interval of time before approaching to the equilibrium which is known as \emph{transient chaos}. When $r \leq 22.35$, all trajectories approach an equilibrium \cite{WJ}.

\begin{figure}[H]
      \centering

      \resizebox{14cm}{!}
      {
      \begin{tabular}{c c c}
      \includegraphics [scale=8]{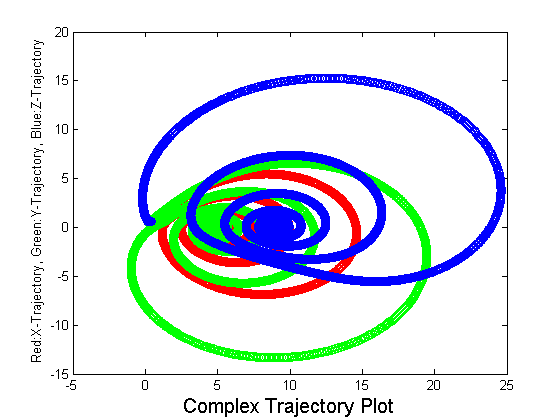}
      \includegraphics [scale=8]{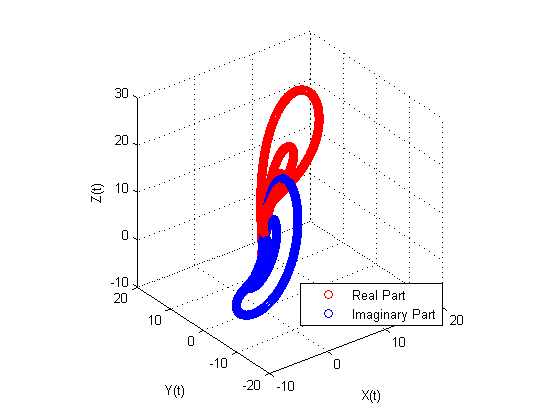}
      \includegraphics [scale=8]{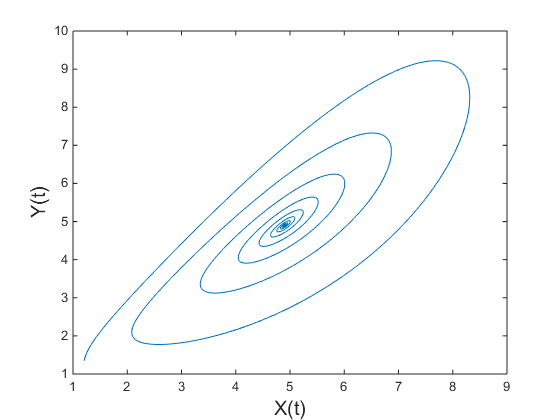}
      \end{tabular}
      }
\caption{For r=10: Trajectory Plots (Left: A complex trajectory (Eq.(3,4,5)), Middle: Trajectory plot in 3 dimension (Eq.(3,4,5)) and Right: Trajectory of the discrete real Lorenz System)}
      \begin{center}

      \end{center}
      \end{figure}
\noindent
In the case, when $r=10$, in both the discrete Lorenz system (Real and Complex), all trajectories are approaching toward one of the equilibriums Converges to one of the sinks $(4.899,4.8990,9)$ which happens to be a sink as stated in Table $1$ and shown in Fig. $3$.\\

\noindent
Considering $r=18$, it is found in the discrete complex Lorenz system Eq. $(3,4,5)$, the trajectories approach one of the equilibriums $(-6.74,-6.74,17)$ without any transient chaos whereas in the discrete real Lorenz system there is transient chaos which is reported in \cite{WJ} and as shown in Fig. $4$.

\begin{figure}[H]
      \centering

      \resizebox{16cm}{!}
      {
      \begin{tabular}{c c c}
      \includegraphics [scale=8]{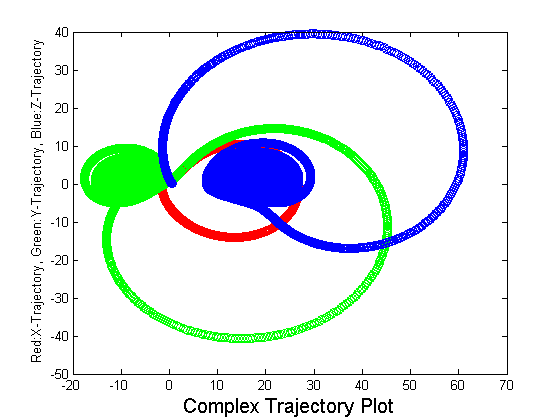}
      \includegraphics [scale=8]{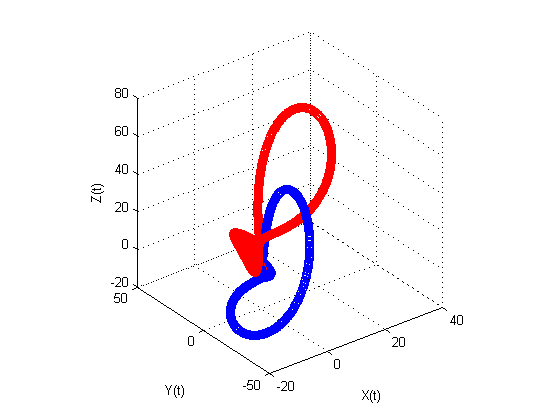}
      \includegraphics [scale=8]{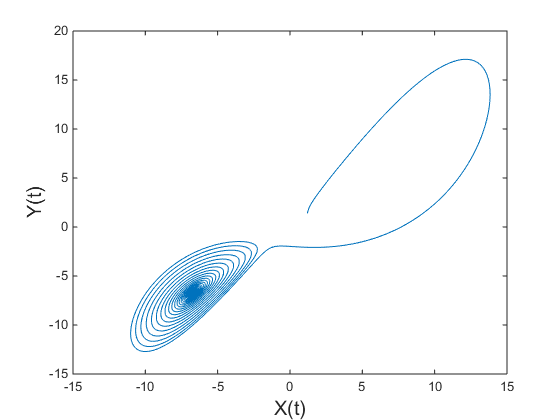}

            \end{tabular}
      }
\caption{For r=18: Trajectory Plots (Left: A complex trajectory (Eq.(3,4,5)), Middle: Trajectory plot (Red: Real Part, Blue: Imaginary Part) in 3 dimension (Eq.(3,4,5)) and Right: Trajectory of the discrete real Lorenz System)}
      \begin{center}

      \end{center}
      \end{figure}

\begin{figure}[H]
      \centering

      \resizebox{15.5cm}{!}
      {
      \begin{tabular}{c c}
      \includegraphics [scale=8]{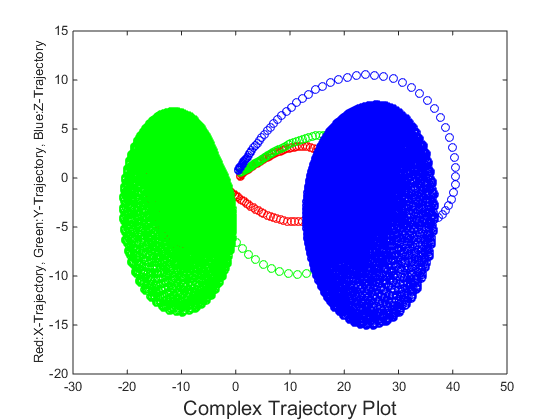}
      \includegraphics [scale=8]{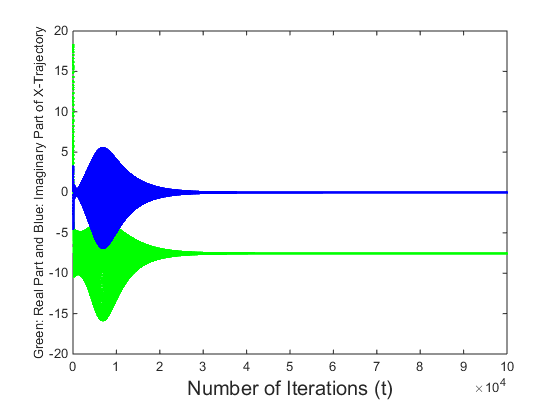}
      \includegraphics [scale=8]{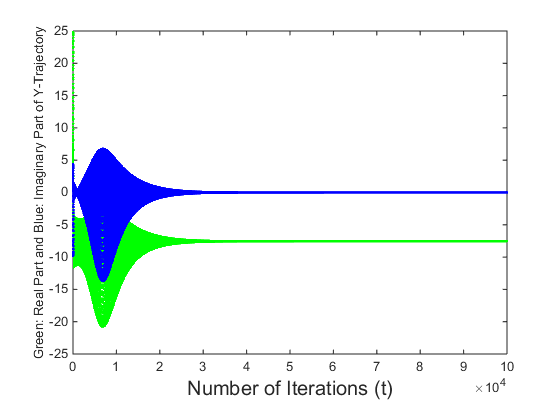}\\
      \includegraphics [scale=8]{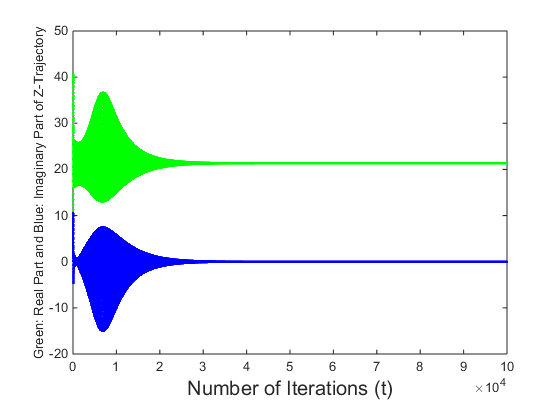}
      \includegraphics [scale=8]{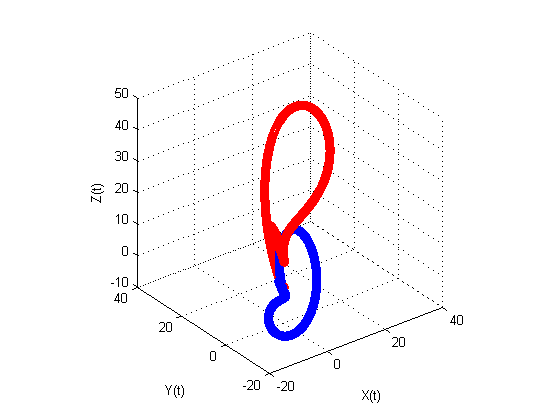}
      \includegraphics [scale=8]{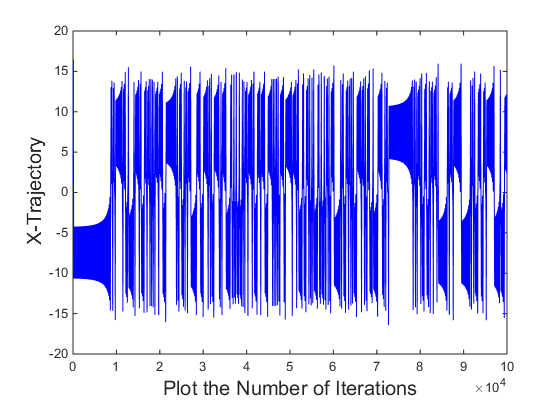}\\
      \includegraphics [scale=8]{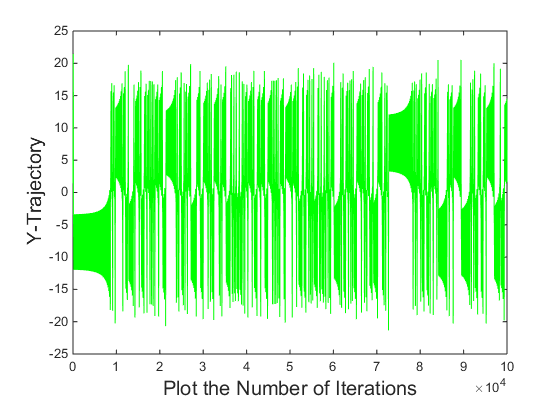}
      \includegraphics [scale=8]{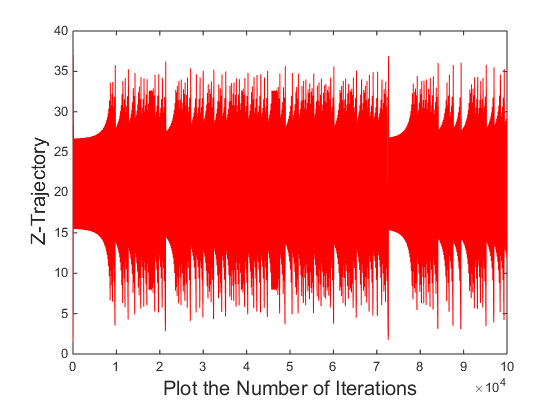}
      \includegraphics [scale=8]{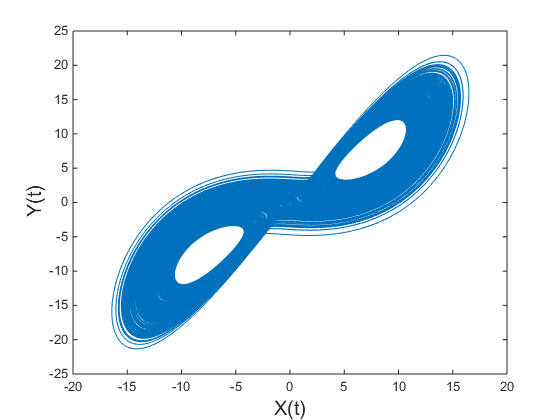}\\

            \end{tabular}
      }
\caption{For r=22.35: Top Left: Complex trajectory plots, Top Middle: X-trajectory, Top Right: Y-trajectory, Middle Left: Z-Trajectory Middle Middle: Trajectory plot in 3 dimension, (In real discrete Lorenz System:: Middle Right: X trajectory, Bottom Left: Y trajectory, Bottom Middle: Z trajectory, Bottom Right: Two dimensional trajectory plot in real discrete system.)}

      \end{figure}

\noindent
 When $r=22.35$, it is seen that in the discrete complex Lorenz system Eq. $(3,4,5)$, the trajectories approach one of the equilibriums $(-7.55, -7.55, 21.35)$ whereas in the discrete real Lorenz system the trajectories possess chaos as shown in Fig.$5$.

\begin{figure}[H]
      \centering

      \resizebox{16.5cm}{!}
      {
      \begin{tabular}{c c}
      \includegraphics [scale=8]{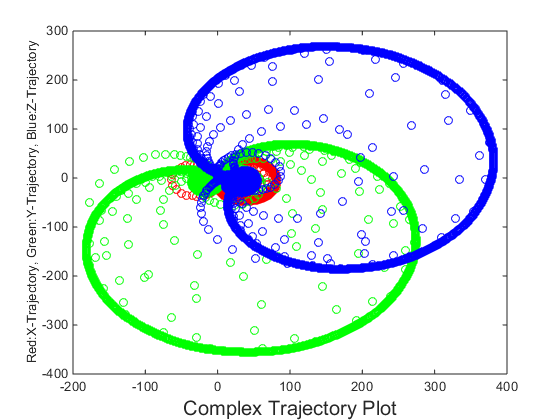}
      \includegraphics [scale=8]{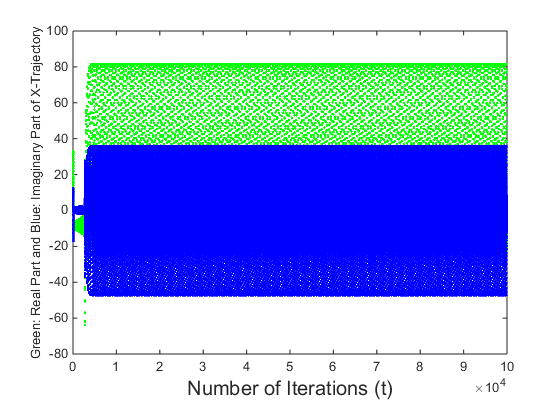}
      \includegraphics [scale=8]{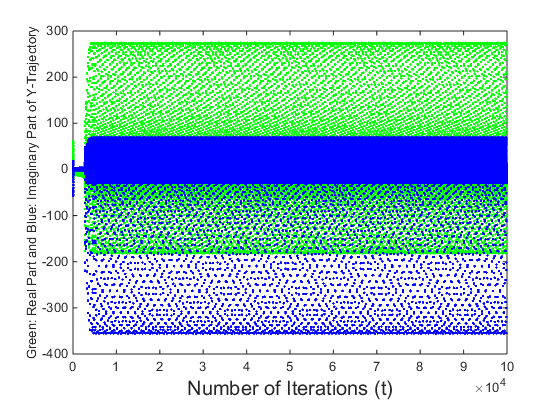}\\
      \includegraphics [scale=8]{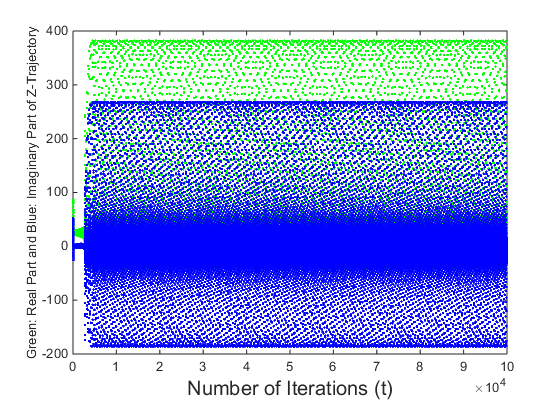}
      \includegraphics [scale=8]{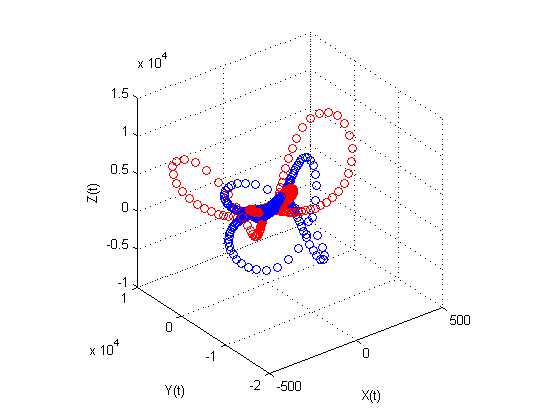}
      \includegraphics [scale=8]{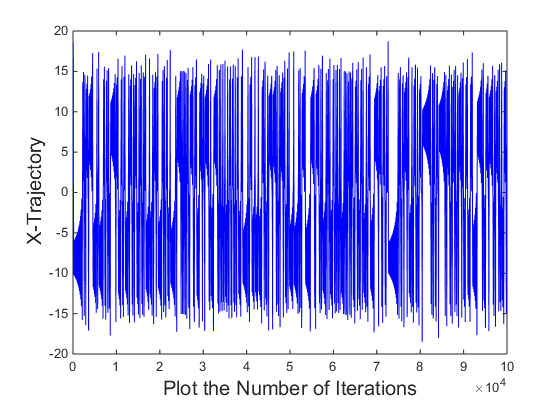}\\
      \includegraphics [scale=8]{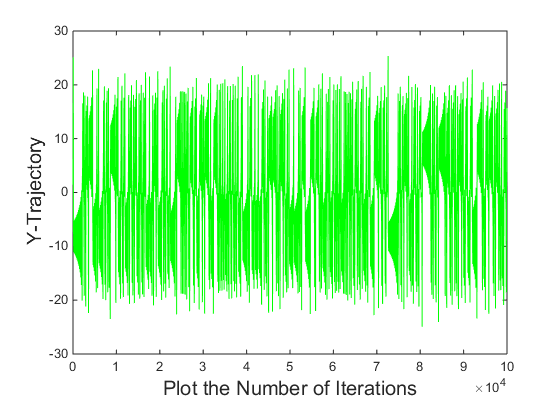}
      \includegraphics [scale=8]{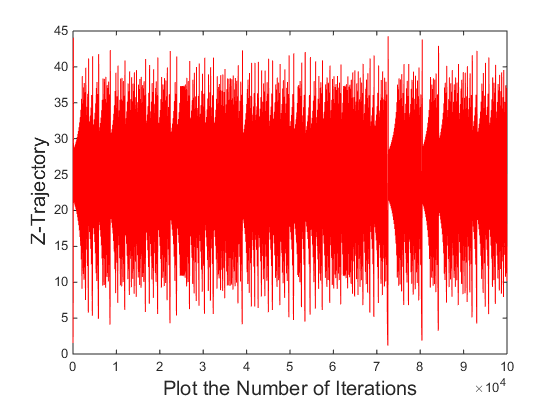}
      \includegraphics [scale=8]{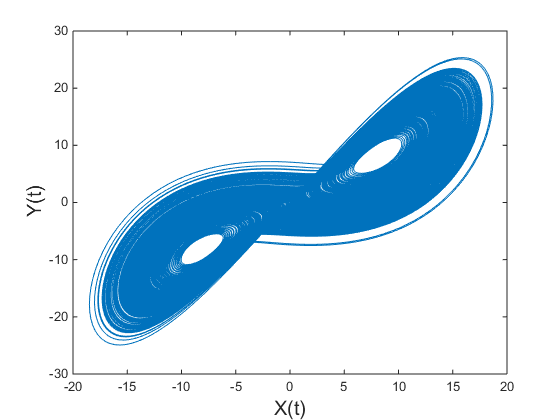}

            \end{tabular}
      }
\caption{For r=26: Top Left: Complex trajectory plots, Top Middle: X-trajectory, Top Right: Y-trajectory, Middle Left: Z-Trajectory Middle Middle: Trajectory plot in 3 dimension, (In Discrete real Lorenz System:: Middle Right: X trajectory, Bottom Left: Y trajectory, Bottom Middle: Z trajectory, Bottom Right: Two dimensional trajectory plot in real discrete system.) }
      \begin{center}

      \end{center}
      \end{figure}
\noindent
When $r=26$, it is observed that in the discrete complex Lorenz system Eq. $(3,4,5)$ and discrete real Lorenz system, the trajectories are chaotic as shown in Fig.$6$.\\

\noindent
 When $r=28$, it is seen that in the discrete complex Lorenz system, the trajectories show transient chaos where as in the discrete real Lorenz system the trajectories possess chaos as shown in Fig.$7$.

\begin{figure}[H]
      \centering

      \resizebox{16.5cm}{!}
      {
      \begin{tabular}{c c c}
      \includegraphics [scale=8]{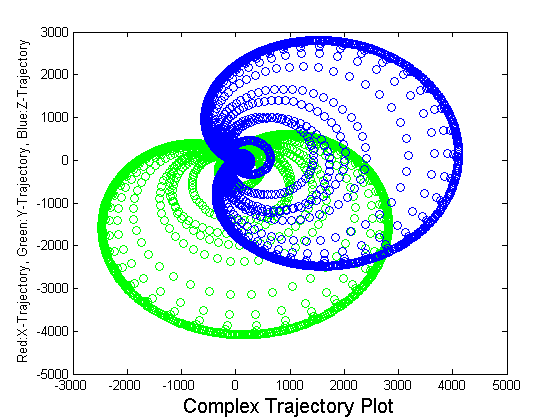}
      \includegraphics [scale=8]{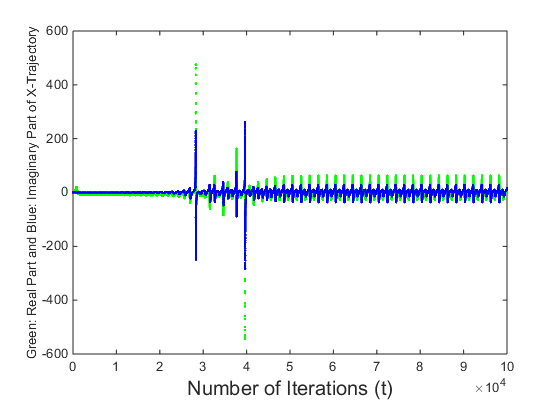}
      \includegraphics [scale=8]{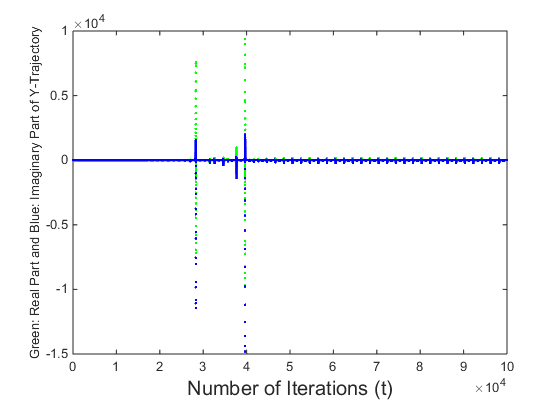}\\
      \includegraphics [scale=8]{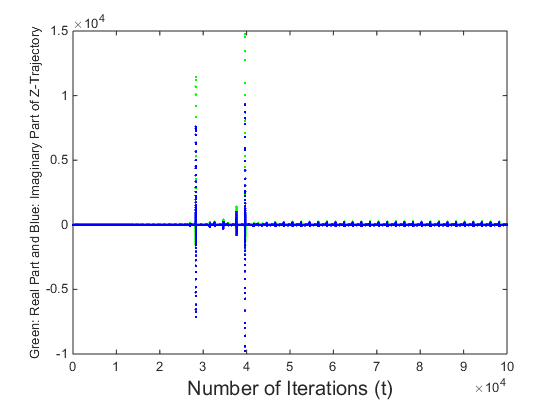}
      \includegraphics [scale=8]{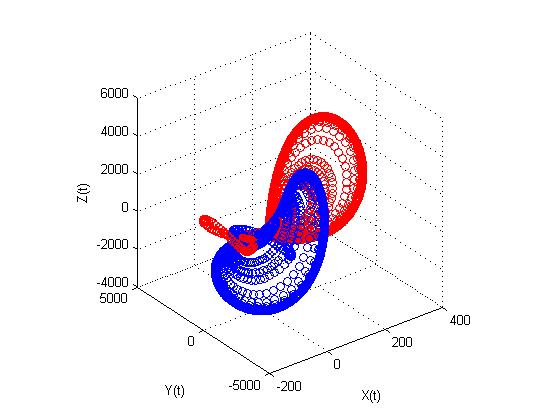}
      \includegraphics [scale=8]{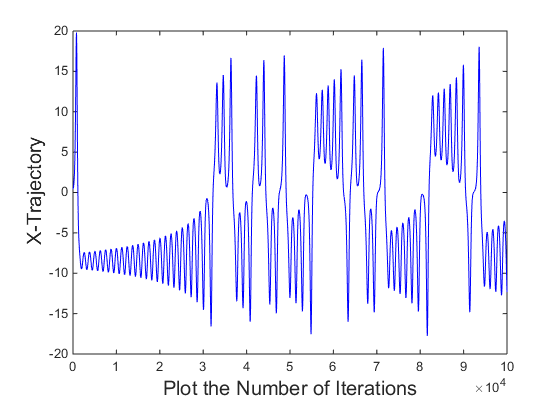}\\
      \includegraphics [scale=8]{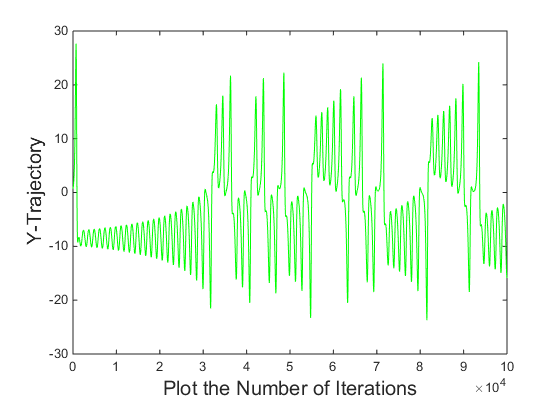}
      \includegraphics [scale=8]{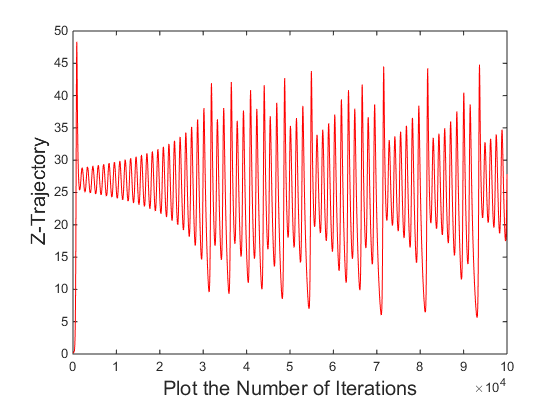}
      \includegraphics [scale=8]{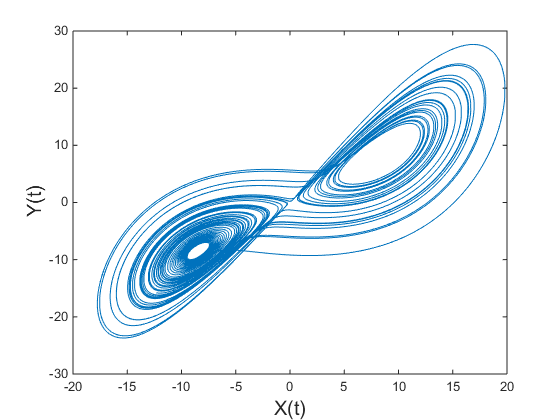}

            \end{tabular}
      }
\caption{For r=28: Top Left: Complex trajectory plots, Top Middle: X-trajectory, Top Right: Y-trajectory, Middle Left: Z-Trajectory Middle Middle: Trajectory plot in 3 dimension, (In Discrete real Lorenz System:: Middle Right: X trajectory, Bottom Left: Y trajectory, Bottom Middle: Z trajectory, Bottom Right: Two dimensional trajectory plot in real discrete system.) }
      \begin{center}

      \end{center}
      \end{figure}
\noindent
In the following subsection we are going explore the dynamics of the complex classical Lorenz system and the complex discrete Lorenz system vividly and compare the dynamics.

\subsection{Comparison of Dynamics of the Complex Discrete Lorenz System and Complex Classical Lorenz System}

Here for different values of the control parameter $r$, the dynamics of complex Lorenz system have been seen and compared with the existing classical complex Lorenz system. Before we proceed further to compare, it is important to note that there is an existing comparison between the classical real Lorenz system and complex Lorenz system and in brief the results are given in the Table $2$ \cite{MH}:\\

\noindent
Here we took the initial value of for both the system $(x_0, y_0, z_0)=(0.1+0.2i, 0.3+0.4i, 1+2i)$ and we fix the parameter a=10, $b=\frac{8}{3}$ and $dt=0.0005$. We took one sample value of control parameter $r$ from each of the specified interval of $r$ in the Table $2$, and run the discrete complex Lorenz system and the result are given in the Table $3$.

\begin{table}[H]

\begin{tabular}{| m{1cm}||m{3.5cm}| |m{4.5cm}| |m{4.5cm}|}
\hline \centering \textbf{Serial No} &
\begin{center}
\textbf{Control Parameter $r$}
\end{center}
 &
\begin{center}
\textbf{Behavior in Real Classical Lorenz System}
\end{center}
&
\begin{center}
\textbf{Behavior in Complex Classical Lorenz System}
\end{center}
\\
\hline \centering 1 &
\centering $1-3.1$
&
\begin{center}
Fixed Point
\end{center}
&
\begin{center}
Periodic
\end{center}\\
\hline
\centering 2 &
\centering $3.1-3.5$
&
\centering  Fixed Point

&
\begin{center}
Quasi Periodic

\end{center}
\\
\hline
\centering 3 &
\centering $3.6-4.2$
&
\begin{center}
 Fixed Point
\end{center}
&
\begin{center}
Chaotic
\end{center}
\\
\hline
\centering 4 &
\centering $4.3-23.9$
&
\centering  Fixed Point

&
\begin{center}
Hyper Chaotic
\end{center}
\\
\hline
\centering 5 &
\centering $23.9-24.7$
&
\begin{center}
Transient Chaotic (Largest Lyapunav exponent: +Ve)
\end{center}
&
\begin{center}
Hyper Chaotic
\end{center}
\\
\hline
\centering 6 &
\centering $24.7-26.7$
&
\begin{center}
One positive Lyapunav exponent
\end{center}
&
\begin{center}
Hyper Chaotic
\end{center}
\\
\hline
\centering 7 &
\centering $26.8-27.5$
&
\begin{center}
One positive Lyapunav exponent
\end{center}
&
\begin{center}
Chaotic
\end{center}
\\
\hline
\centering 8 &
\centering $27.6-27.8$
&
\begin{center}
One positive Lyapunav exponent
\end{center}
&
\begin{center}
Periodic
\end{center}
\\
\hline
\centering 9 &
\centering $27.9-99.6$
&
\begin{center}
Chaotic
\end{center}
&
\begin{center}
Hyper Chaotic
\end{center}
\\
\hline
\centering 10 &
\centering $99.6-100$
&
\begin{center}
Periodic
\end{center}
&
\begin{center}
Hyper Chaotic
\end{center}
\\
\hline

\end{tabular}
\caption{Dynamics of complex classical Lorenz System and real classical Lorenz system for different values of $r$.}
\label{Table:1}
\end{table}

\begin{table}[H]

\begin{tabular}{| m{1cm}||m{3.5cm}| |m{4.5cm}| |m{4.5cm}|}
\hline \centering \textbf{Serial No} &
\begin{center}
\textbf{Control Parameter $r$}
\end{center}
 &
\begin{center}
\textbf{Behavior in Discrete Complex Lorenz System}
\end{center}
&
\begin{center}
\textbf{Behavior in Complex Classical Lorenz System}
\end{center}
\\
\hline \centering 1 &
\centering $1$
&
\centering  Converges to one of the sinks $(0.0554, 0.0554, 0.0012)$

&
\begin{center}
Periodic

\end{center}\\
\hline
\centering 2 &
\centering $2$
&
\centering  Converges to one of the sinks $(1.633, 1.633, 1)$

&
\begin{center}
Periodic

\end{center}
\\
\hline
\centering 3 &
\centering $3$
&
\begin{center}
 Converges to one of the sinks $(2.3094, 2.3094, 2)$
\end{center}
&
\begin{center}
Quasi-periodic
\end{center}
\\
\hline
\centering 4 &
\centering $3.2$
&
\centering   Converges to one of the sinks $(2.4221, 2.4221, 2.2)$

&
\begin{center}
Quasi-periodic
\end{center}
\\
\hline
\centering 5 &
\centering $3.4$
&
\centering  Converges to one of the sinks $(2.5298, 2.5298, 2.4)$

&
\begin{center}
Quasi-periodic
\end{center}
\\
\hline
\centering 6 &
\centering $3.7$
&
\centering  Converges to one of the sinks $(2.6833, 2.6833, 2.7)$

&
\begin{center}
Chaotic
\end{center}
\\
\hline
\centering 7 &
\centering $4.1$
&
\centering  Converges to one of the sinks $(2.8752, 2.8752, 3.1)$

&
\begin{center}
Chaotic
\end{center}
\\
\hline
\centering 8 &
\centering $20$
&
\centering  Converges to one of the sinks $(-7.1181, -7.1181, 19)$

&
\begin{center}
Hyper-Chaotic
\end{center}
\\
\hline
\centering 9 &
\centering $24$
&
\centering  Transient Chaos (One positive Lyapunav exponent)

&
\begin{center}
Hyper-Chaotic
\end{center}
\\
\hline
\centering 10 &
\centering $24.9$
&
\centering  Transient Chaos (One positive Lyapunav exponent)

&
\begin{center}
Hyper-Chaotic
\end{center}
\\
\hline
\centering 11 &
\centering $26$
&
\centering  Transient Chaos (One positive Lyapunav exponent)

&
\begin{center}
Hyper-Chaotic
\end{center}
\\
\hline
\centering 12 &
\centering $26.7$
&
\centering  Chaos (Largest Lyapunav exponent is positive)

&
\begin{center}
Hyper-Chaotic
\end{center}
\\
\hline
\centering 13 &
\centering $27$
&
\centering  Chaos (Largest Lyapunav exponent is positive)

&
\begin{center}
Chaos
\end{center}
\\
\hline
\centering 14 &
\centering $27.7$
&
\centering  Chaos (Largest Lyapunav exponent is positive)

&
\begin{center}
Periodic
\end{center}
\\
\hline
\centering 15 &
\centering $100$
&
\centering  Transient Chaos (One positive Lyapunav exponent)

&
\begin{center}
Hyper-Chaos
\end{center}
\\
\hline

\end{tabular}
\caption{Dynamics of complex Lorenz System and complex classical Lorenz System for different values of $r$.}
\label{Table:2}
\end{table}

\noindent

\begin{figure}[H]
      \centering

      \resizebox{10cm}{!}
      {
      \begin{tabular}{c c c}
      \includegraphics [scale=8]{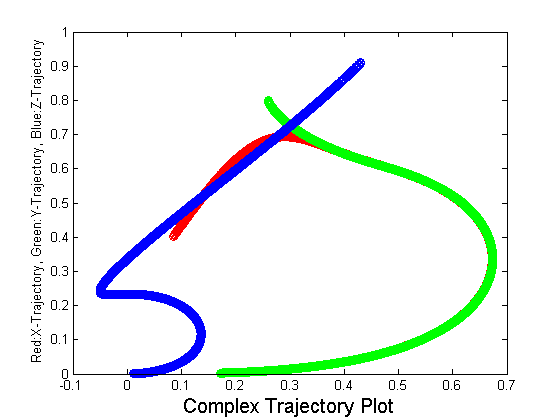}
      \includegraphics [scale=8]{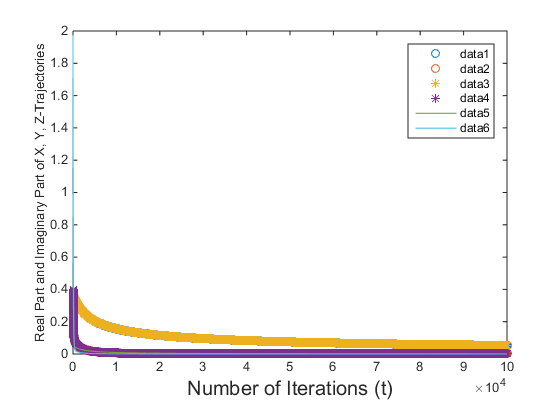}
      \includegraphics [scale=8]{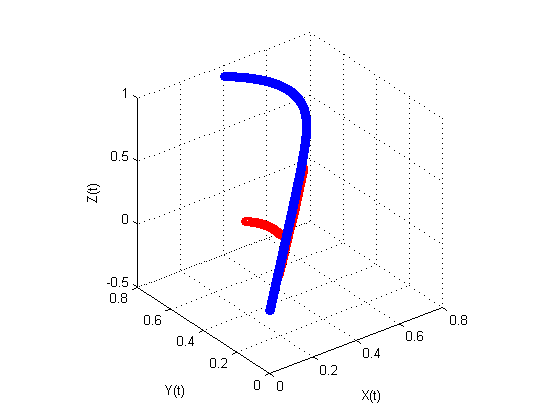}\\
      \includegraphics [scale=8]{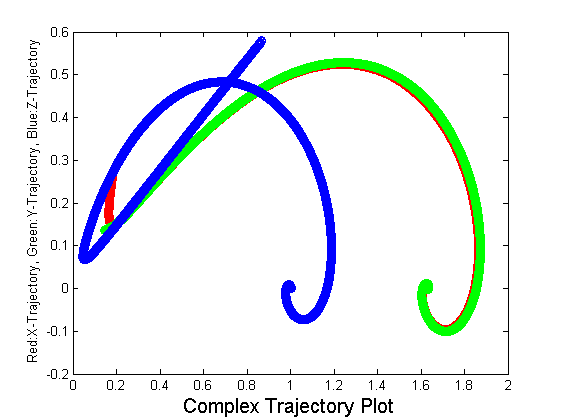}
      \includegraphics [scale=8]{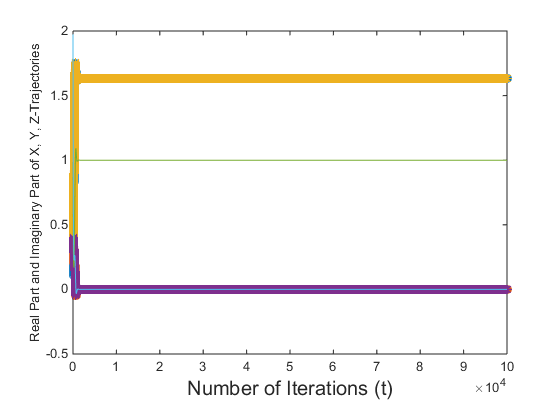}
      \includegraphics [scale=8]{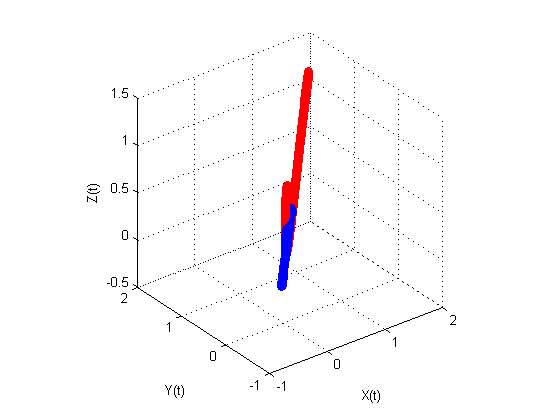}\\
      \includegraphics [scale=8]{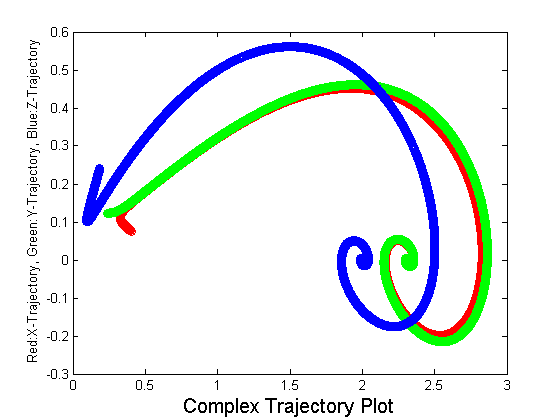}
      \includegraphics [scale=8]{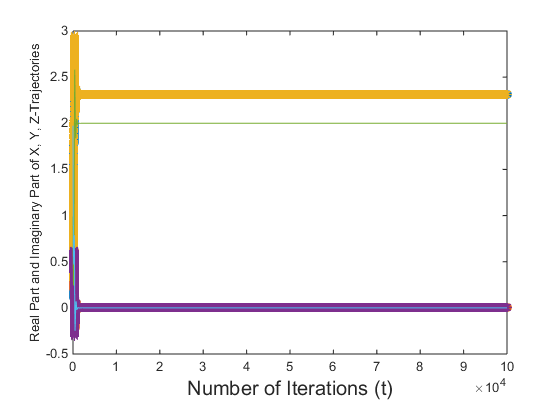}
      \includegraphics [scale=8]{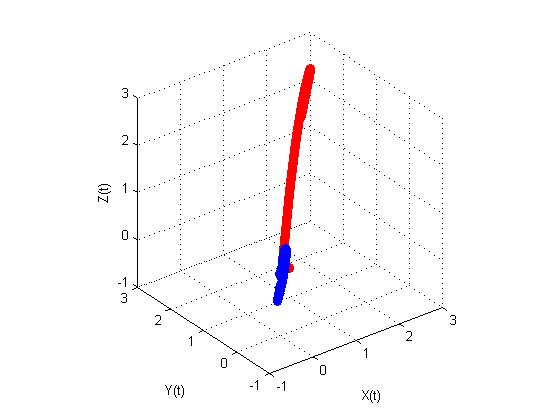}\\
      \includegraphics [scale=8]{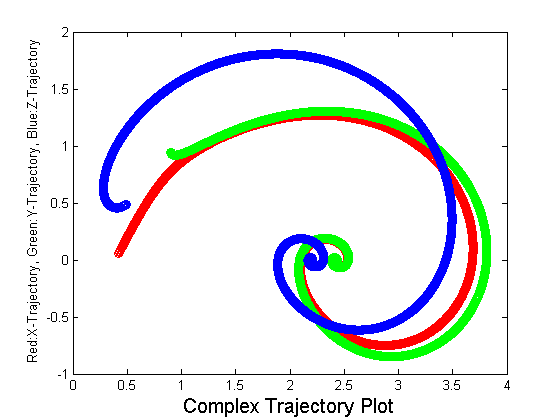}
      \includegraphics [scale=8]{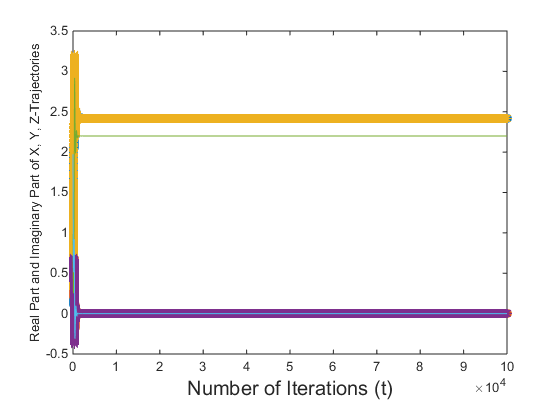}
      \includegraphics [scale=8]{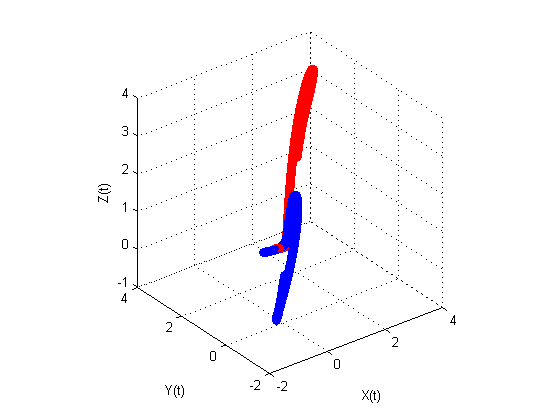}\\
      \includegraphics [scale=8]{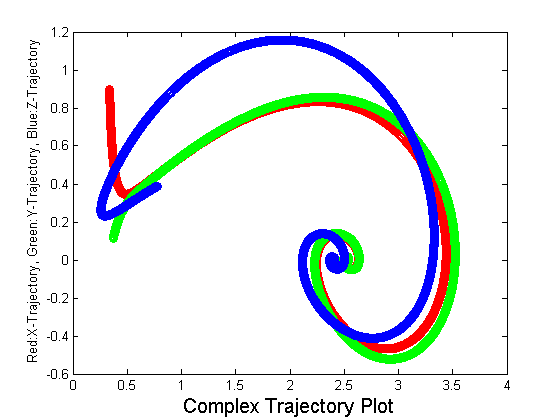}
      \includegraphics [scale=8]{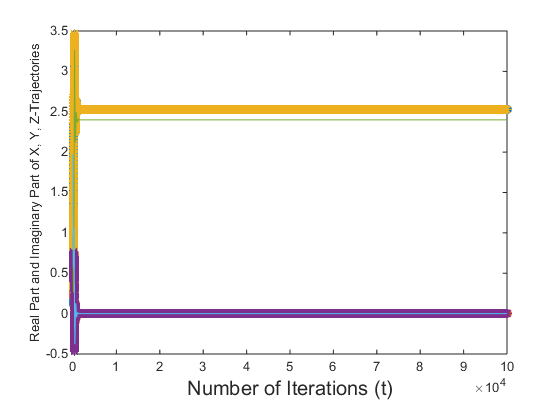}
      \includegraphics [scale=8]{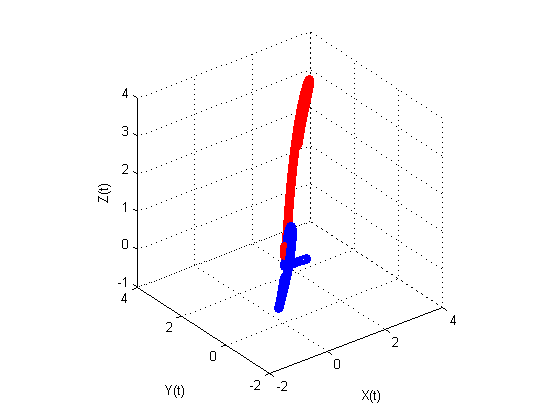}\\
      \includegraphics [scale=8]{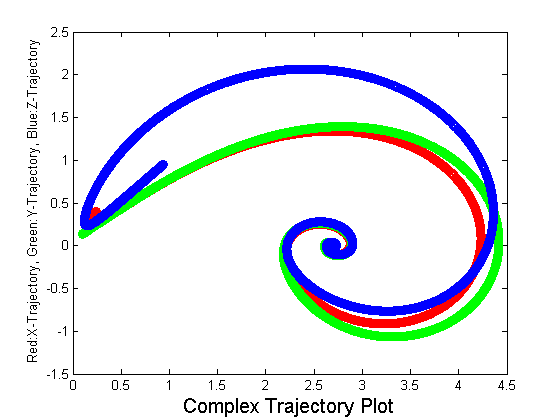}
      \includegraphics [scale=8]{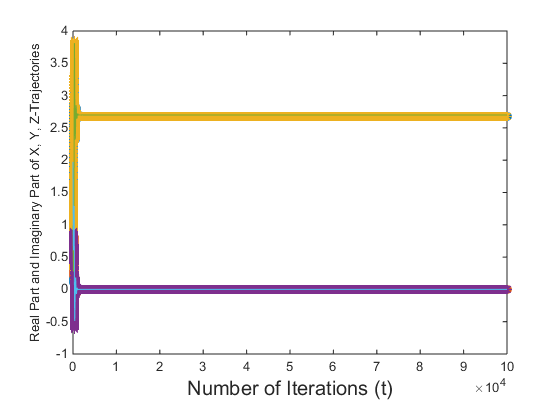}
      \includegraphics [scale=8]{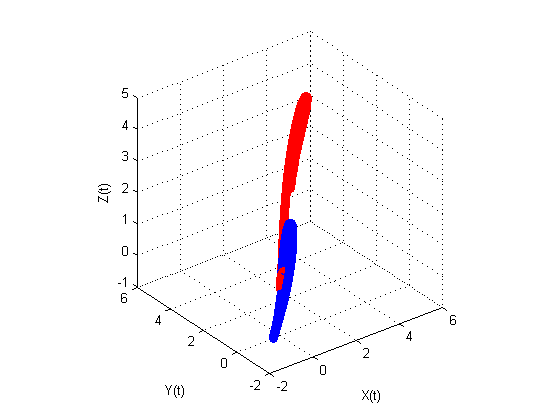}\\
      \includegraphics [scale=8]{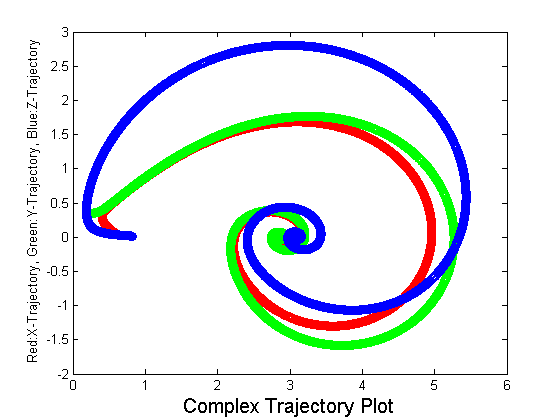}
      \includegraphics [scale=8]{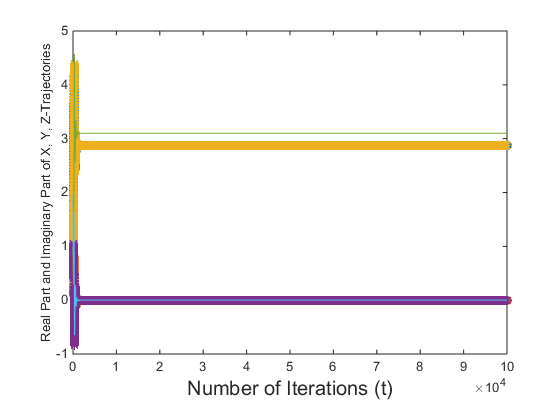}
      \includegraphics [scale=8]{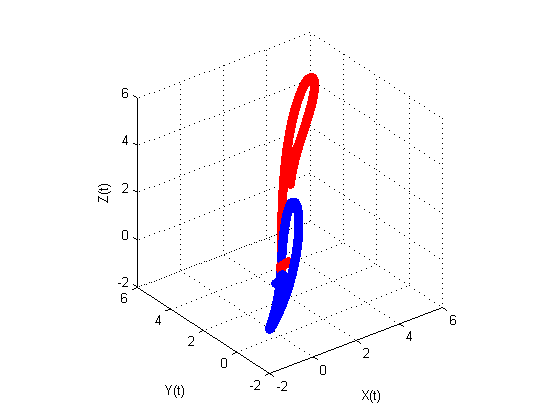}\\
      \includegraphics [scale=8]{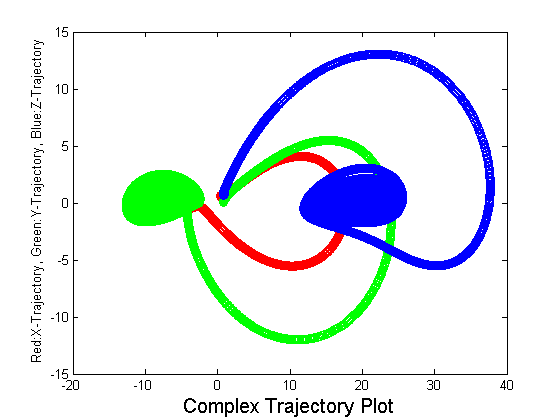}
      \includegraphics [scale=8]{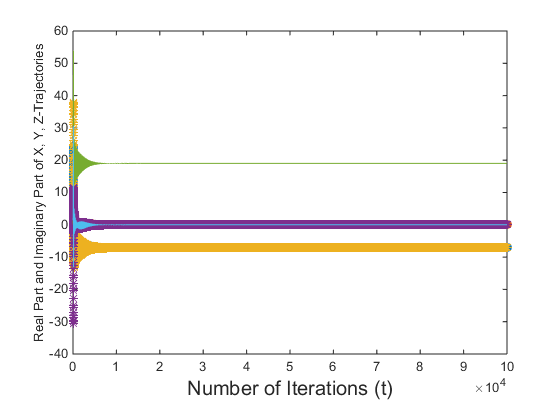}
      \includegraphics [scale=8]{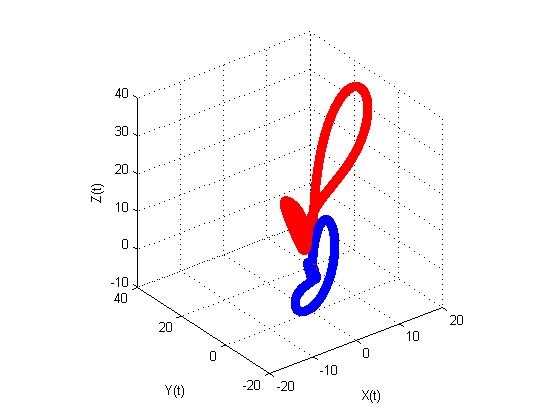}\\

            \end{tabular}
      }
\caption{For r=1, 2, 3, 3,2, 3.4, 3.7, 4.1, 20 respectively from top to bottom row: Top Left: Complex trajectory plots, Top Middle: X, Y, Z-trajectories, Top Right: Trajectory plot (Red: Real part, Blue: Imaginary part) in 3 dimension. }
      \begin{center}

      \end{center}
      \end{figure}

\begin{figure}[H]
      \centering

      \resizebox{12cm}{!}
      {
      \begin{tabular}{c c c}
      \includegraphics [scale=8]{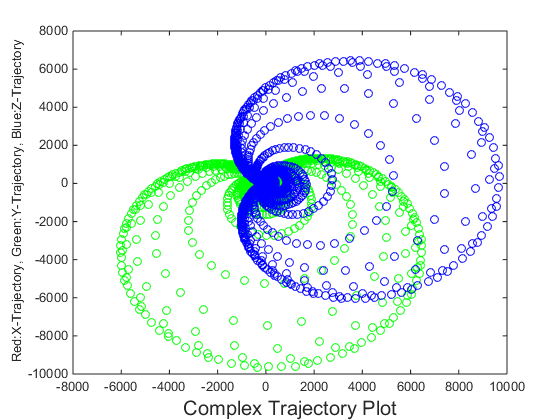}
      \includegraphics [scale=8]{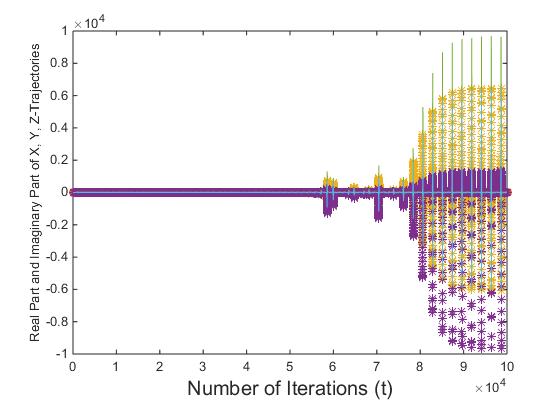}
      \includegraphics [scale=8]{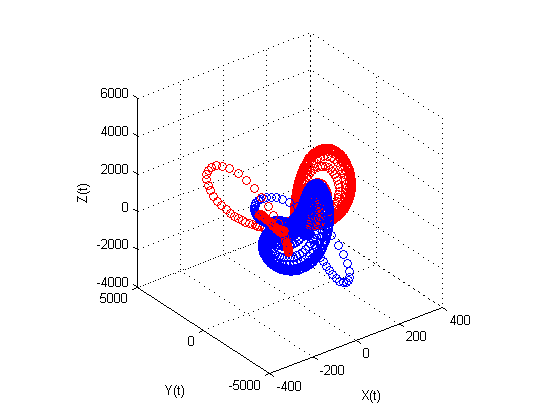}\\
      \includegraphics [scale=8]{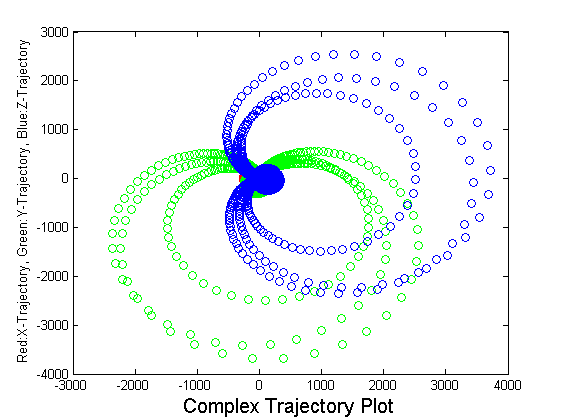}
      \includegraphics [scale=8]{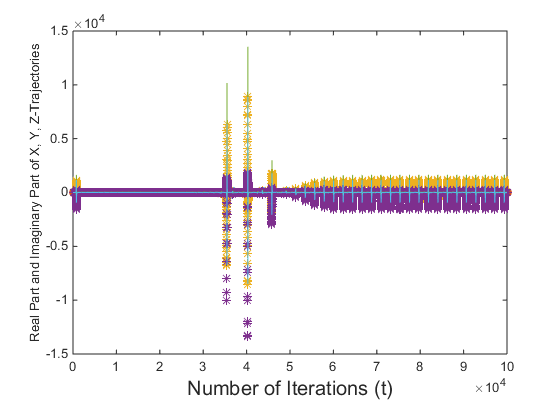}
      \includegraphics [scale=8]{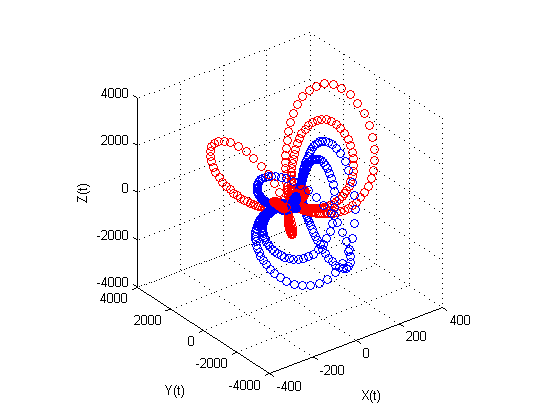}\\
      \includegraphics [scale=8]{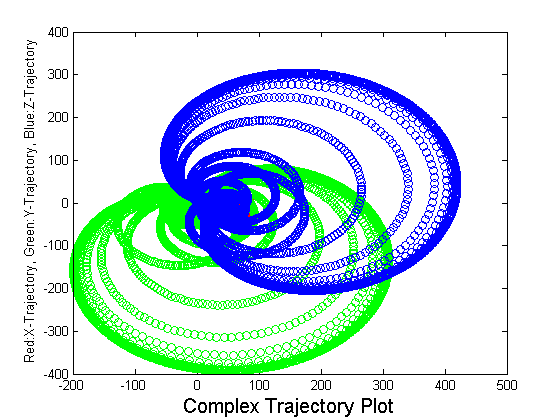}
      \includegraphics [scale=8]{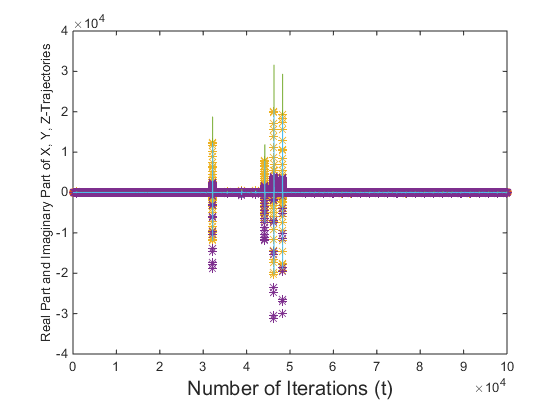}
      \includegraphics [scale=8]{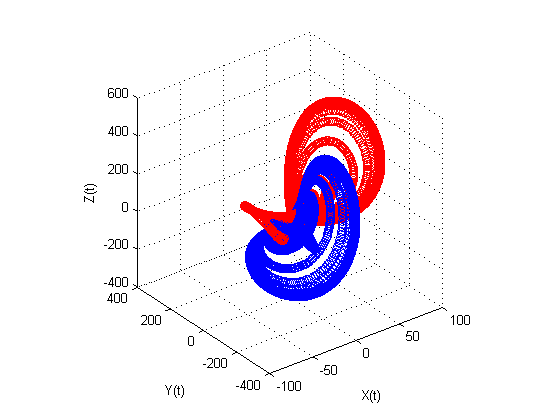}\\
      \includegraphics [scale=8]{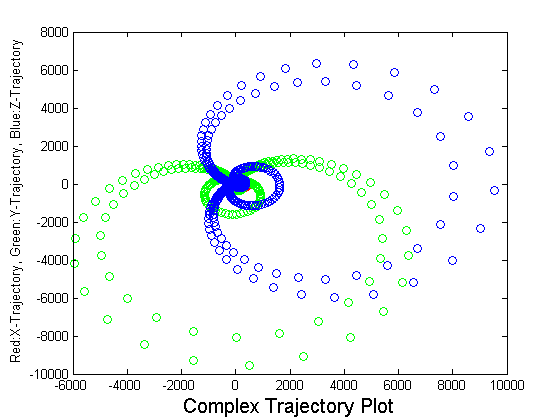}
      \includegraphics [scale=8]{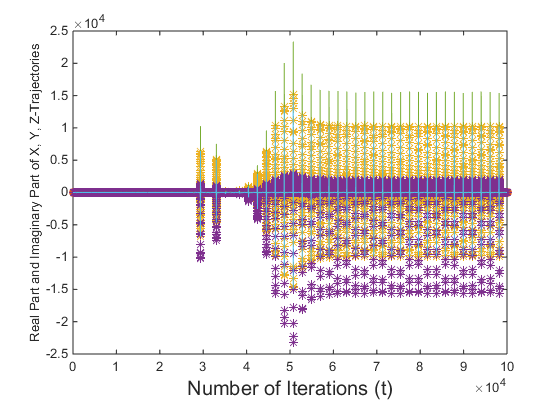}
      \includegraphics [scale=8]{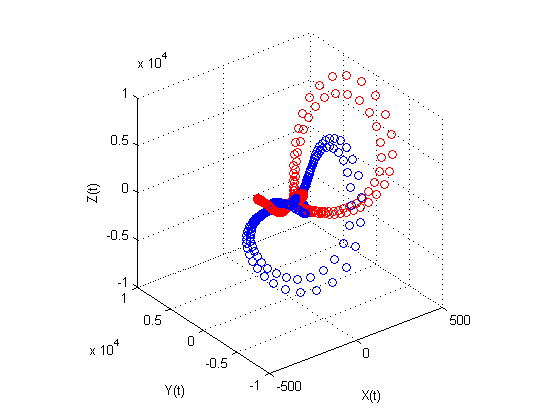}\\
      \includegraphics [scale=8]{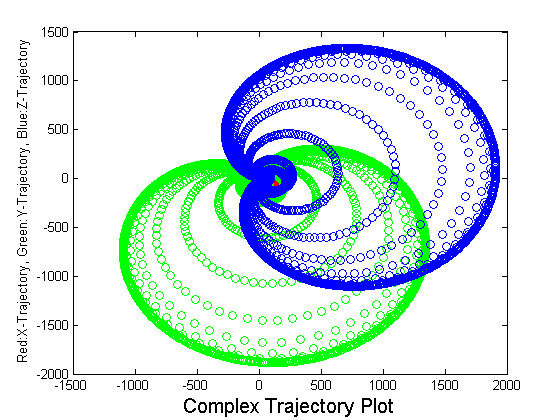}
      \includegraphics [scale=8]{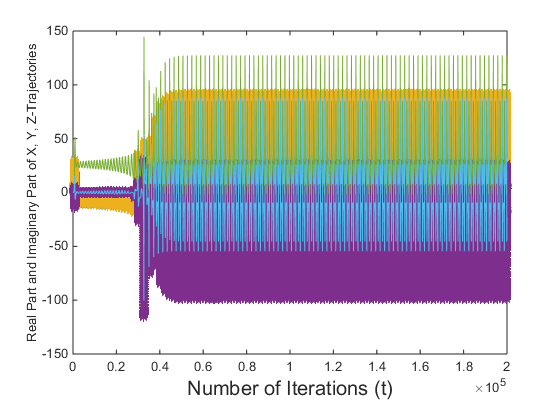}
      \includegraphics [scale=8]{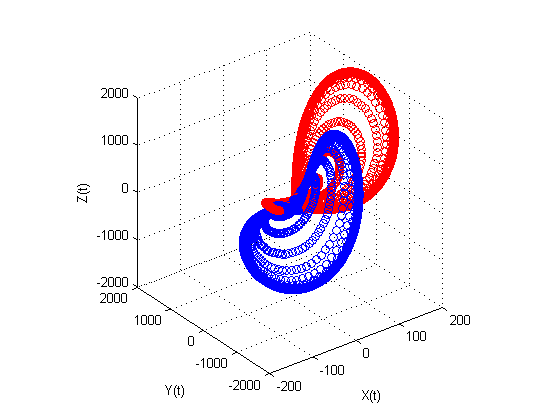}\\
      \includegraphics [scale=8]{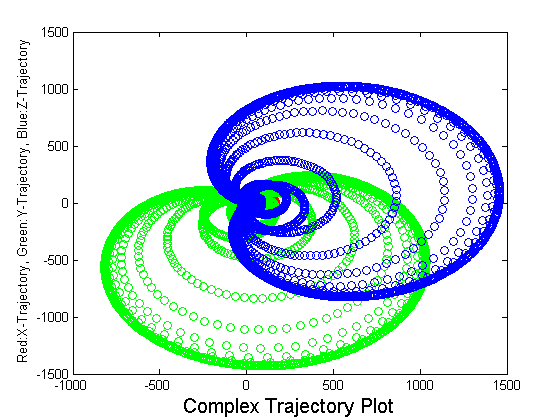}
      \includegraphics [scale=8]{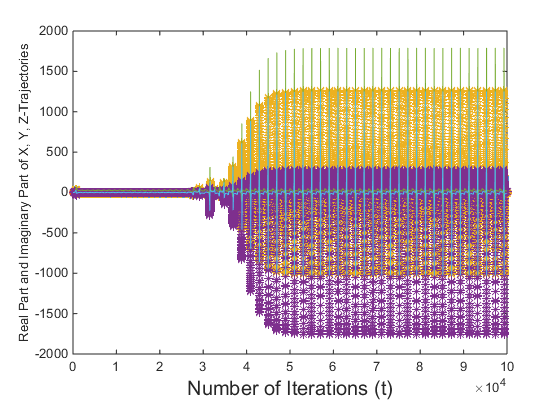}
      \includegraphics [scale=8]{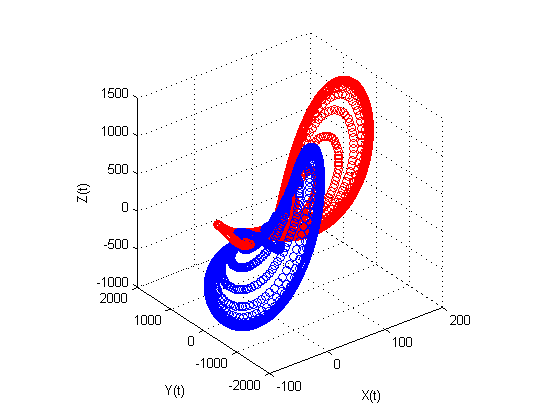}\\
      \includegraphics [scale=8]{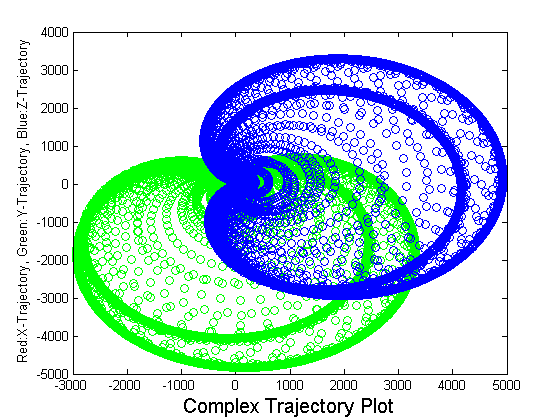}
      \includegraphics [scale=8]{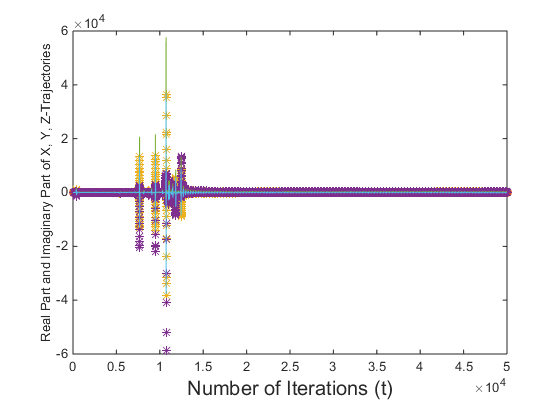}
      \includegraphics [scale=8]{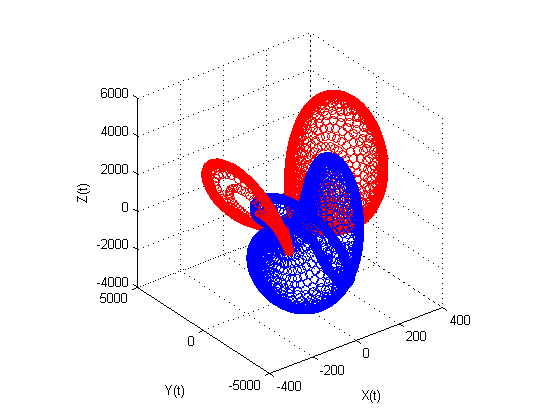}\\

            \end{tabular}
      }
\caption{For r=24, 24.9, 26, 26.7, 27, 27.7, 100 respectively from top to bottom row: Top Left: Complex trajectory plots, Top Middle: X, Y, Z-trajectories, Top Right: Trajectory plot (Red: Real part, Blue: Imaginary part) in 3 dimension. }
      \begin{center}

      \end{center}
      \end{figure}

\section{Complex Control Parameters}
Here we are about to explore the dynamics of the discrete complex Lorenz system Eq. (3,4,5) considering first the control parameter as origin (zero) in the complex plane and in second assuming the parameter $r$ is a non-origin which are depicted as the following.

\subsection{Case when the Parameter $r=0$}

First we consider $r=0$ and $dt=0.0005$ and see the dynamics as depicted in the following Table $4$ and Fig.$10$. When $r=0$, the fixed points of the Lorenz system Eq. (3,4,5) become $(0,0,0), (-i\sqrt{b}, -i\sqrt{b},-1)$ and $(i\sqrt{b}, i\sqrt{b},-1)$.\\

\noindent
In the Table $4$, the serial number $1$, the parameters are setting as $r=0$, $dt=0.0005$, $a=4+20i$ and $b=20+10i$, the fixed points $(0,0,0)$ and $(\pm (4.60221 + 1.08643i), \pm (4.60221 + 1.08643i), 1)$. The real part of the eigenvalues of the jacobian about the fixed points are all positive and therefore the three fixed points are sink (attracting). Under the initial condition the trajectory converges to the sink $(0,0,0)$ which is shown in the Table $4$ and Fig. $10$.

\begin{table}[H]

\begin{tabular}{| m{1cm}||m{3.8cm}| |m{4.5cm}| |m{4.5cm}|}
\hline \centering \textbf{Serial No} &
\begin{center}
\textbf{Parameter $a$ and $b$}
\end{center}
 &
\begin{center}
\textbf{Behavior in Discrete Complex Lorenz System}
\end{center}
&
\begin{center}
\textbf{Dynamics in the Classical Real Lorenz System}
\end{center}
\\
\hline
\centering 1 &
\centering $a=4+20i$, $b=20+10i$
&
\centering  Converges to (0,0,0)

&
\begin{center}
Either unbounded or two attractors are coexisting

\end{center}
\\
\hline
\centering 2 &
\centering $a=0.052$, $b=-0.3$, Note: (-1,0,1) is the initial value taken.
&
\begin{center}
Triangle like attractor with very high periodicity.
\end{center}
&
\begin{center}
Two attractor coexist.
\end{center}
\\
\hline
\centering 3 &
\centering $a=0.256$, $b=-0.3$
&
\begin{center}
Two attractors coexist.
\end{center}
&
\begin{center}
Two attractors coexist.
\end{center}
\\
\hline
\centering 4 &
\centering $a=0.277$, $b=-0.3$
&
\begin{center}
Very high period
\end{center}
&
\begin{center}
Two attractors coexist.
\end{center}
\\
\hline

\end{tabular}
\caption{Dynamics of Discrete Complex Lorenz System and Classical Real Lorenz System when $r=0$.}
\label{Table:1}
\end{table}

\noindent
With the setting $r=0$, $dt=0.0005$, $a=0.052$ and $b=-0.3$ as in the serial number $2$ of the Table $4$, the fixed points $(0,0,0)$ and $(\pm 0.547723i, \pm 0.547723i, 1)$. The real part of the eigenvalues of the jacobian about the fixed points $(0,0,0)$ and $(\pm 0.547723i, \pm 0.547723i, 1)$ are all positive and therefore the three fixed points are sink (attracting). But surprisingly, for the said set of parameters and the initial value $(-1,0,1)$, the trajectory converging to a very high periodic triangle-like trajectory which is been depicted in the Table $4$ and Fig. $10$.

\noindent
As we set in serial number $3$ in the Table $4$, the fixed points $(0,0,0)$ and $(\pm 0.547723i, \pm 0.547723i, 1)$ which is same as in the case in serial no $2$. The real part of the eigenvalues of the jacobian about the fixed points $(0,0,0)$ and $(\pm 0.547723i, \pm 0.547723i, 1)$ are all positive and therefore the three fixed points are sink (attracting). Surprisingly, with initial value $(-0.1,+0.1,-2)$, the trajectory converging to a double chaotic attractors which do coexist which is been depicted in the Fig. $10$. This phenomena is exactly seen in the classical real Lorenz system \cite{MH}.

\noindent
Here we differ from the previous case just by changing the parameter $a$ as $0.277$ and as the parameter $a$ does not have any role in the fixed points, so the fixed points remain $(0,0,0)$ and $(\pm 0.547723i, \pm 0.547723i, 1)$ and it is found that the real part of the eigenvalues of the jacobian about the fixed points $(0,0,0)$ and $(\pm 0.547723i, \pm 0.547723i, 1)$ are all positive and naturally so the three fixed points are sink (attracting). But with the initial parameter $(0.1,-0.1,-13)$ the trajectory possesses to a very high periodic attractor as shown in Fig. $10$.

\begin{figure}[H]
      \centering

      \resizebox{16cm}{!}
      {
      \begin{tabular}{c c c c}
      \includegraphics [scale=8]{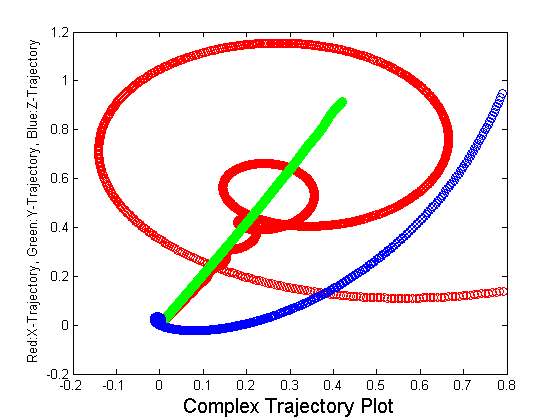}
      \includegraphics [scale=8]{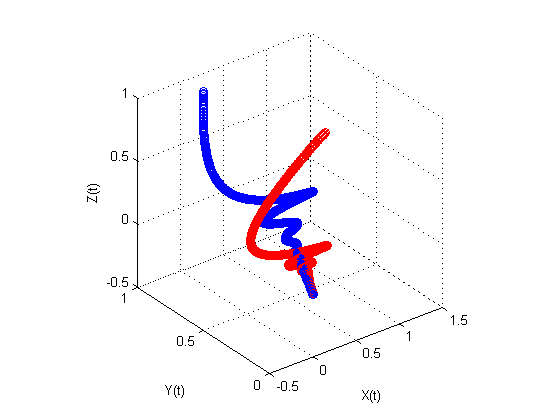}
      \includegraphics [scale=8]{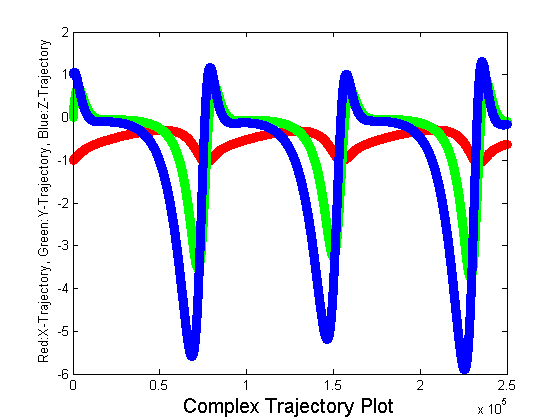}
      \includegraphics [scale=8]{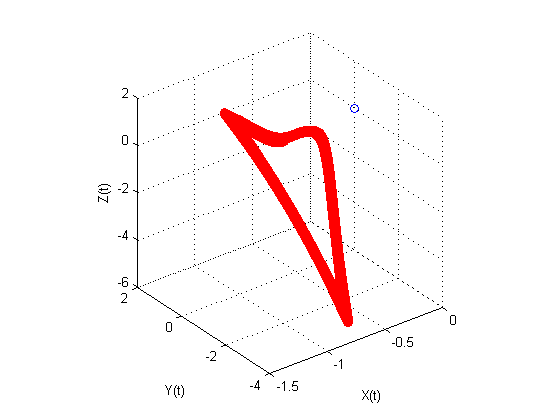}\\
      \includegraphics [scale=8]{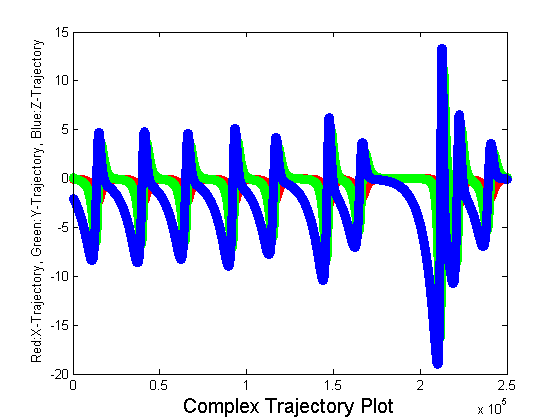}
      \includegraphics [scale=8]{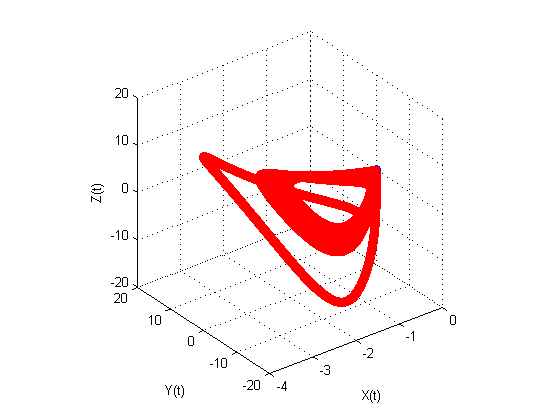}
      \includegraphics [scale=8]{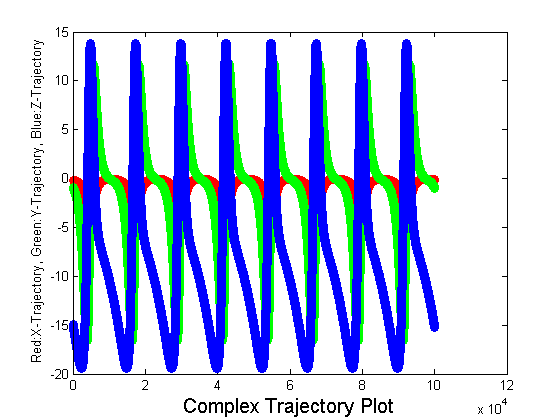}
      \includegraphics [scale=8]{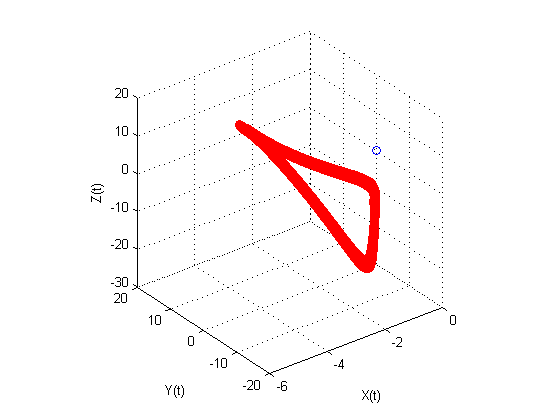}\\

            \end{tabular}
      }
\caption{For r=0: Top left two figures: Case in Serial No.1,Top right two figures: Case in Serial No.2, Bottom left two figures: Case in Serial No.3, Bottom right two figures: Case in Serial No.4, (Red: Real part, Blue: Imaginary part) in 3 dimension. }
    \end{figure}

\noindent
It is noted that, when the parameters $r$ is zero and $b$ is a negative real number there exists a range of values of $a$ for which the trajectory of the classical real Lorenz system proceed to a double coexisting attractor which is matched in the discrete complex Lorenz system only in the case of serial no 2 as stated in the Table $4$ \cite{JCS7}.

\subsection{Case when the Parameter $r$ is a Non-zero}

\noindent
When the control parameter $r$ is a non-zero complex number, we shall see the dynamics of the discrete complex Lorenz system Eq. $(3,4,5)$. \\\\
\noindent
For the parameters as sited in the serial number $1, 2, 4, 6, 7$ and $8$ in the Table 5, we found that $\abs{r}>\abs{a\frac{a+b+3}{a-b-1}}$, this violates the condition for convergence to the non-zero fixed points. It is found that the under these set of parameters the trajectories are converging eventually to the zero fixed point whereas in some of these cases transient chaos have been observed as shown in the Fig. 11 and Fig. 12.\\\\
\noindent
Considering the parameters in the serial number $3, 5, 9$ and $10$, it is found that $\abs{r}<\abs{a\frac{a+b+3}{a-b-1}}$, this ensures the stability of the non-zero fixed points. The trajectory in these cases possesses to transient chaos except the last case as shown in serial number 10, the trajectory is gradually diverging as shown in Fig.$12$.

\begin{table}[H]

\begin{tabular}{| m{1cm}||m{3.8cm}| |m{4.5cm}| |m{4.5cm}|}
\hline \centering \textbf{Serial No} &
\begin{center}
\textbf{Parameters: $a$, $b$ and $r$}
\end{center}
 &
\begin{center}
\textbf{Dynamics in the Discrete Complex Lorenz System}
\end{center}
&
\begin{center}
\textbf{Remark}
\end{center}
\\
\hline
\centering 1 &
\centering $a=0.0462 + 0.0971i$, $b=0.6948 + 0.3171i$ and $r=0+9i$
&
\centering Converges to (0,0,0)

&
\begin{center}
The fixed point (0,0,0) is attracting. The real part of the eigenvalues of the jacobian at (0,0,0) are all positive.
\end{center}
\\
\hline
\centering 2 &
\centering $a= 0.0759 + 0.0540i$, $b= 0.9340 + 0.1299i$, and $r=-5 - 8i$
&
\begin{center}
Three attractors coexist.
\end{center}
&
\begin{center}
All three fixed points are attracting. The real part of the eigenvalues of the jacobian at all three fixed points are all positive.
\end{center}
\\
\hline
\centering 3 &
\centering $a=0.5497 + 0.9172i$, $b=0.7572 + 0.7537i
$ and $r=0-4i$
&
\begin{center}
Transient chaos.
\end{center}
&
\begin{center}
Two attractors coexist.
\end{center}
\\
\hline
\centering 4 &
\centering $a=0.6020+0.2630i$, $b=0.6892+0.7482i$ and $r=0+3i$
&
\begin{center}
Converges to (0,0,0).
\end{center}
&
\begin{center}
The fixed point (0,0,0) is attracting. The real part of the eigenvalues of the jacobian at (0,0,0) are all positive.
\end{center}
\\
\hline
\centering 5 &
\centering $a=0.1869 + 0.4898i$, $b=0.6463 + 0.7094i$ and $r=0-i$
&
\begin{center}
Transient Chaos
\end{center}
&
\begin{center}
All the three fixed points are attracting.
\end{center}
\\
\hline
\centering 6 &
\centering $a=0.4733 + 0.3517i$, $b=0.5497 + 0.9172i$ and $r=7+2i$
&
\begin{center}
Transient Chaos
\end{center}
&
\begin{center}
All the three fixed points are attracting.
\end{center}
\\
\hline
\centering 7 &
\centering $a=0.8258 + 0.5383i$, $b=0.4427 + 0.1067i$ and $r=10-9i$
&
\begin{center}
Chaos (Two attractors coexist.)
\end{center}
&
\begin{center}
All the three fixed points are attracting.
\end{center}
\\
\hline
\centering 8 &
\centering $a=0.3998 + 0.2599i$, $b=0.9106 + 0.1818i$ and $r=6-i$
&
\begin{center}
Converges to $(2.1760 - 0.0003i, 2.1760 - 0.0003i, 5-i)$
\end{center}
&
\begin{center}
All the three fixed points are attracting.
\end{center}
\\
\hline
\centering 9 &
\centering $a=0.7702 + 0.3225i$, $b=0.0358 + 0.1759i$ and $r=8+5i$
&
\begin{center}
Chaos (Double attractor coexists.)
\end{center}
&
\begin{center}
All the three fixed points are attracting.
\end{center}
\\
\hline
\centering 10 &
\centering $a=10$, $b=\frac{8}{3}$ and $r=-1-5i$
&
\begin{center}
Gradually divergent
\end{center}
&
\begin{center}
All the three fixed points are attracting.
\end{center}
\\
\hline

\end{tabular}
\caption{Dynamics of Discrete Complex Lorenz System when $r$ is varying over y-axis and punctured complex plane.}
\label{Table:1}
\end{table}

\begin{figure}[H]
      \centering

      \resizebox{16cm}{!}
      {
      \begin{tabular}{c c c}
      \includegraphics [scale=8]{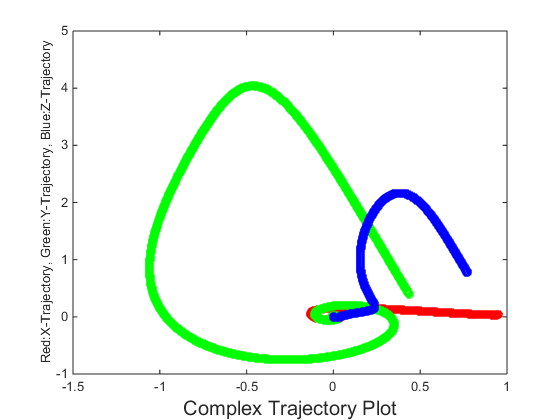}
      \includegraphics [scale=8]{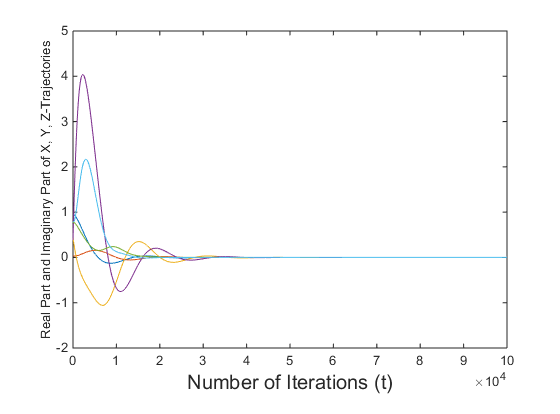}
      \includegraphics [scale=8]{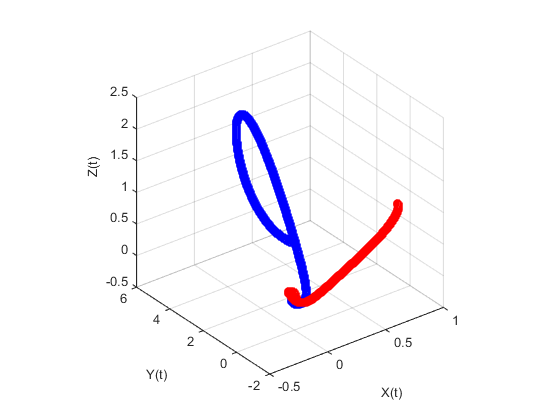}\\
      \includegraphics [scale=8]{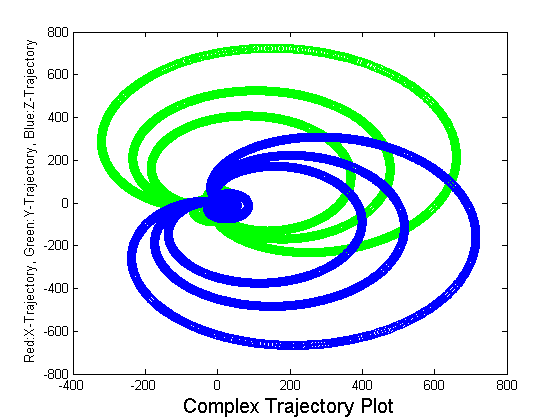}
      \includegraphics [scale=8]{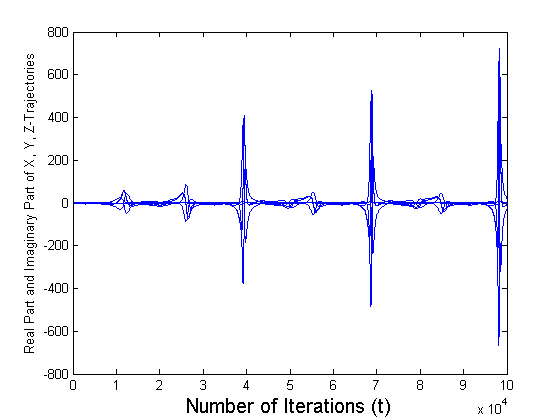}
      \includegraphics [scale=8]{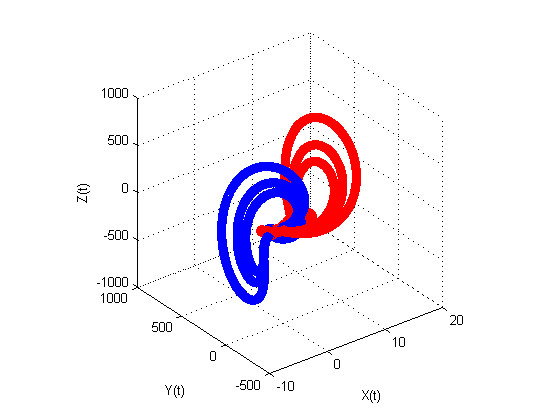}\\
      \includegraphics [scale=8]{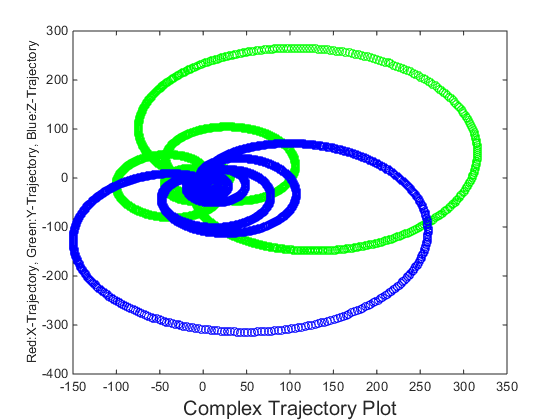}
      \includegraphics [scale=8]{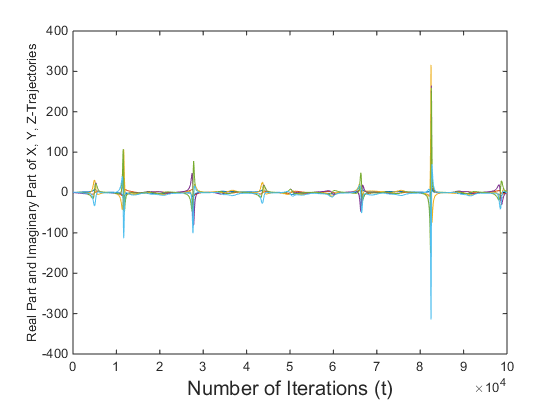}
      \includegraphics [scale=8]{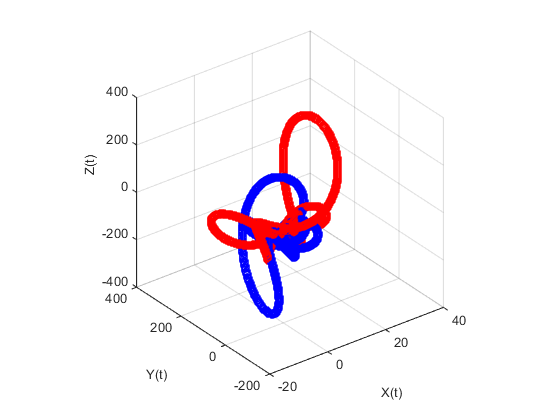}\\
      \includegraphics [scale=8]{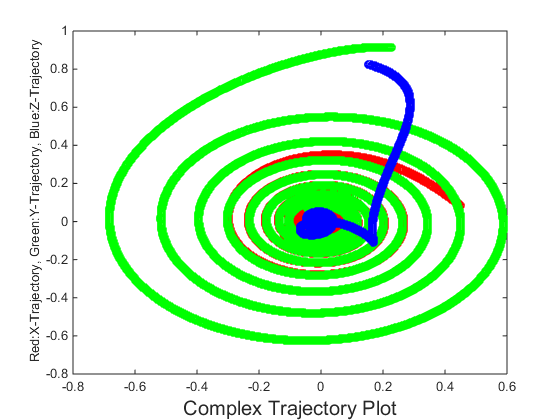}
      \includegraphics [scale=8]{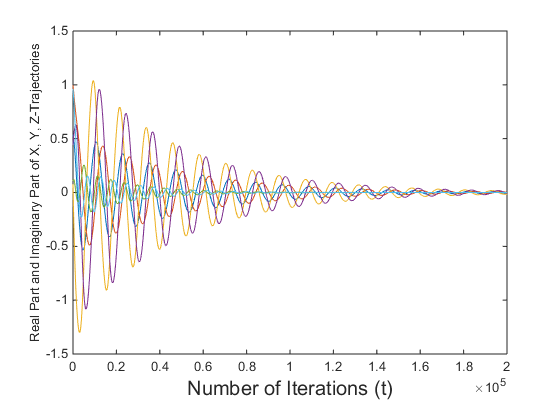}
      \includegraphics [scale=8]{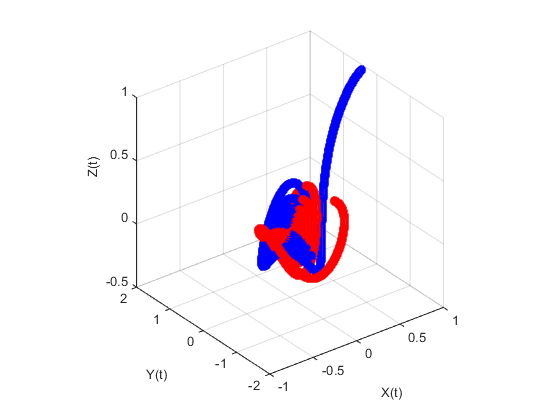}\\
      \includegraphics [scale=8]{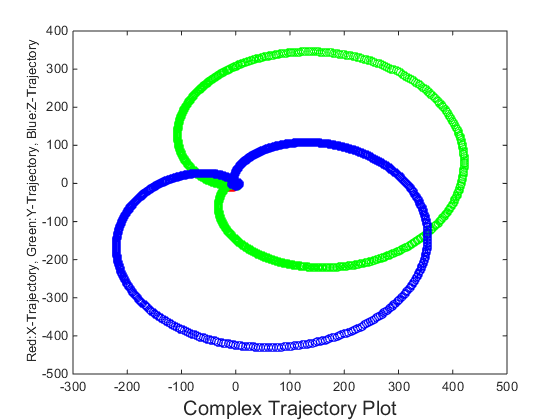}
      \includegraphics [scale=8]{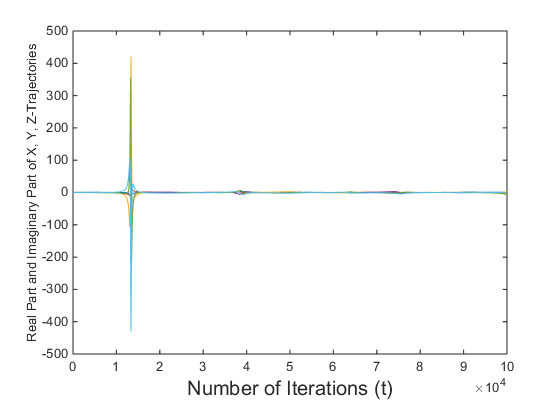}
      \includegraphics [scale=8]{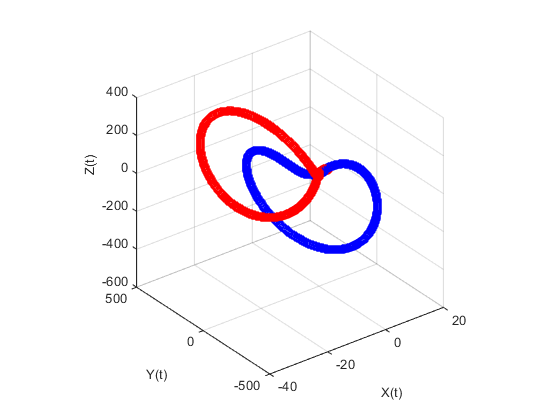}\\

            \end{tabular}
      }
\caption{Serial number $1$ to $5$ in the Table 5, from top to bottom row: Top Left: Complex trajectory plots, Top Middle: X, Y, Z-trajectories, Top Right: Trajectory plot (Red: Real part, Blue: Imaginary part) in 3 dimension.}
      \begin{center}

      \end{center}
      \end{figure}

\begin{figure}[H]
      \centering

      \resizebox{16cm}{!}
      {
      \begin{tabular}{c c c}
      \includegraphics [scale=8]{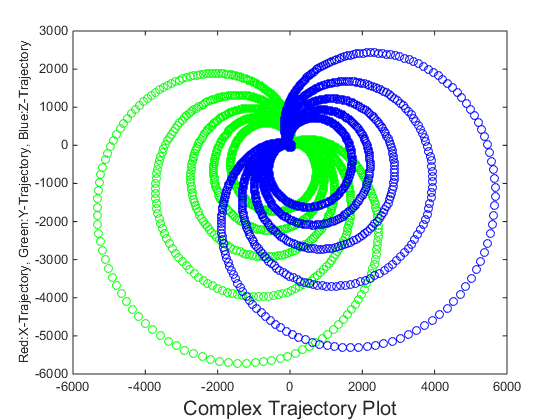}
      \includegraphics [scale=8]{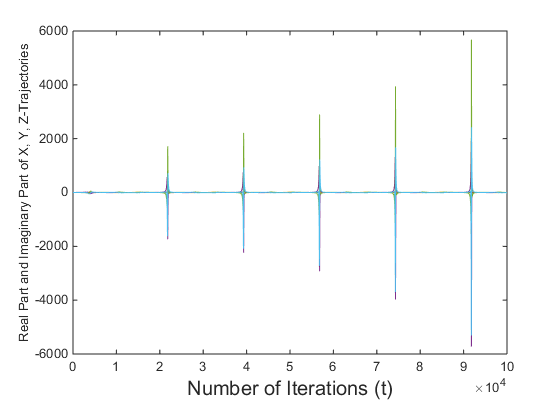}
      \includegraphics [scale=8]{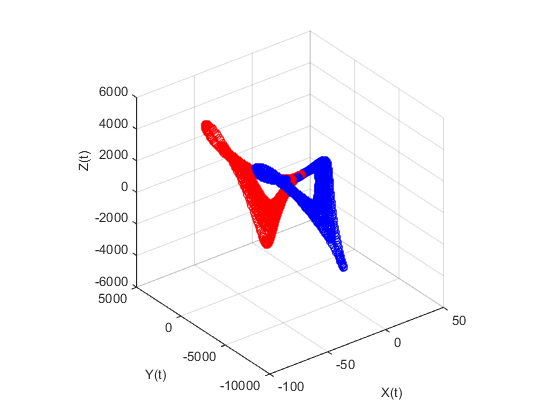}\\
      \includegraphics [scale=8]{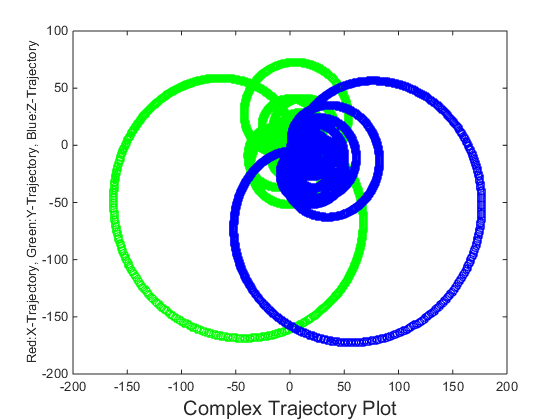}
      \includegraphics [scale=8]{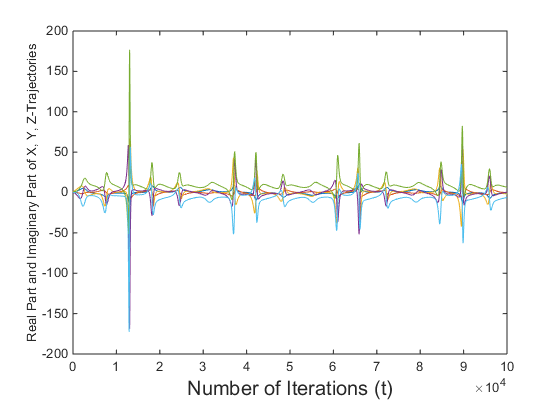}
      \includegraphics [scale=8]{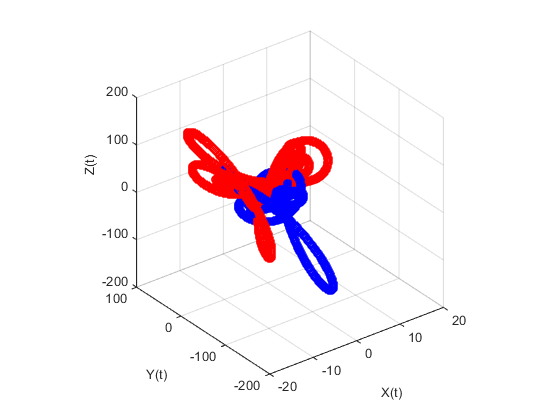}\\
      \includegraphics [scale=8]{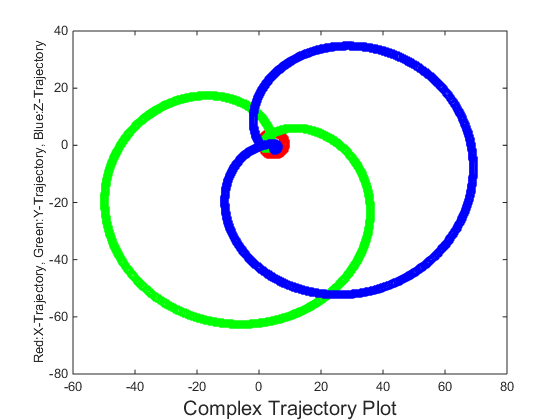}
      \includegraphics [scale=8]{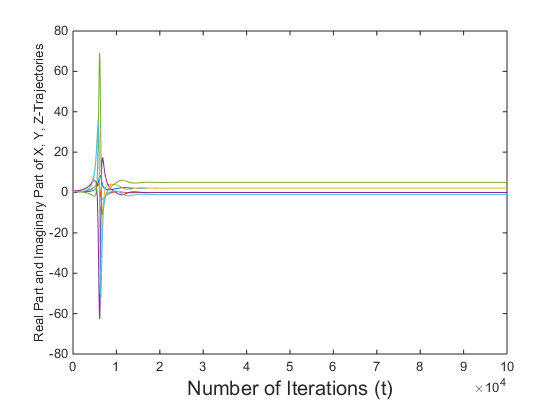}
      \includegraphics [scale=8]{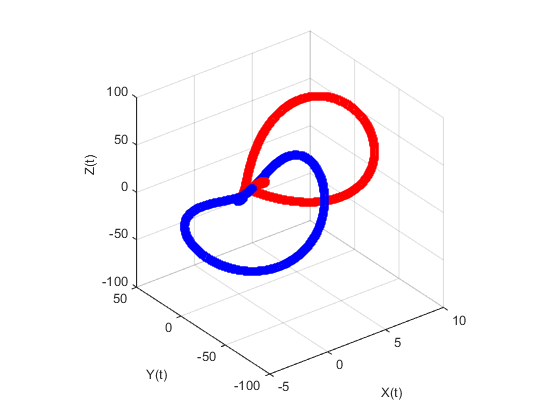}\\
      \includegraphics [scale=8]{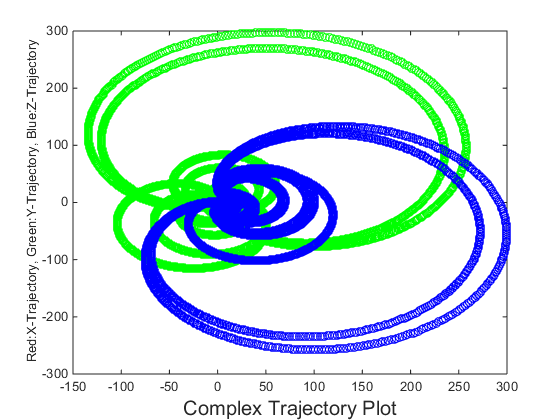}
      \includegraphics [scale=8]{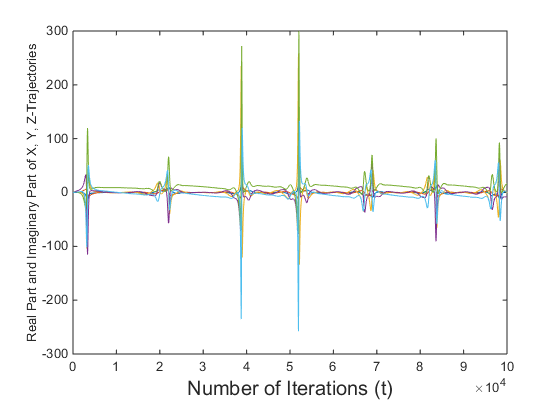}
      \includegraphics [scale=8]{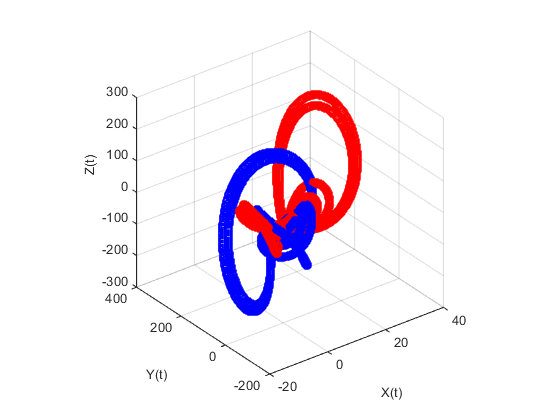}\\
     \includegraphics [scale=8]{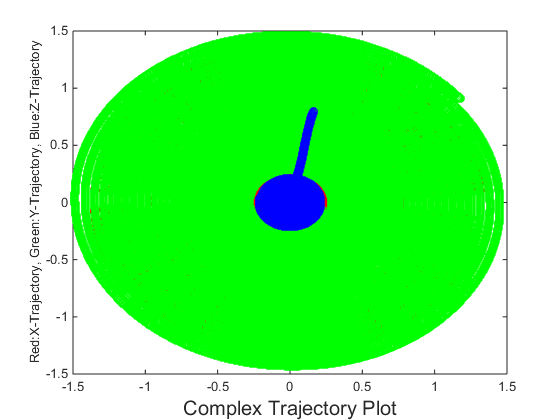}
     \includegraphics [scale=8]{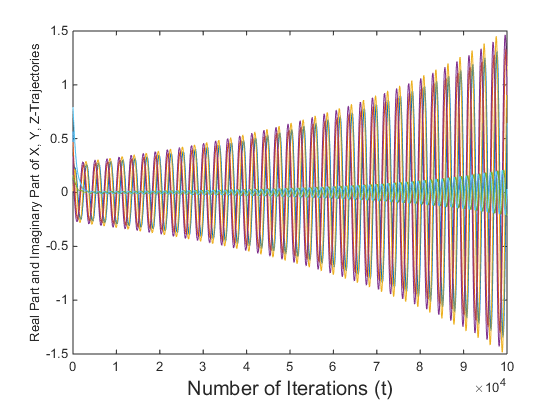}
     \includegraphics [scale=8]{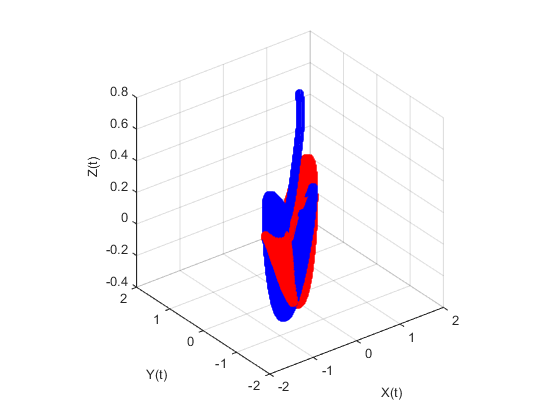}\\

            \end{tabular}
      }
\caption{Serial number $6$ to $10$ in the Table 5, from top to bottom row: Top Left: Complex trajectory plots, Top Middle: X, Y, Z-trajectories, Top Right: Trajectory plot (Red: Real part, Blue: Imaginary part) in 3 dimension. }
      \begin{center}

      \end{center}
      \end{figure}

\section{A Few Special Examples}

In this section, we would like to encounter a set of appealing examples in the discrete complex Lorenz system. It is found in real classical Lorenz system that for $r=99.96$ the system possesses knotted periodic orbits \cite{JCS2}. Here let us see a few examples of dynamics of the Lorenz system Eq. $(3,4,5)$ where $\abs{r}=99.96$ with $dt=0.00005$.

\begin{table}[H]

\begin{tabular}{| m{1cm}||m{3.8cm}| |m{4.5cm}| |m{4.5cm}|}
\hline \centering \textbf{Serial No} &
\begin{center}
\textbf{Parameters: $a$, $b$ and $r$}
\end{center}
 &
\begin{center}
\textbf{Initial Value}
\end{center}
&
\begin{center}
\textbf{Dynamics in the Discrete Complex Lorenz System}
\end{center}
\\
\hline
\centering 1 &
\centering $a= 0.1576 + 0.9706i$, $b=0.9572 + 0.4854i$ and $r=-71.5539 -69.8000i$
&
 \begin{center}
 $( 0.8003 + 0.1419i, 0.4218 + 0.9157i,  0.7922 + 0.9595i)$
 \end{center}
&
\begin{center}
Knotted Chaotic orbit
\end{center}
\\
\hline
\centering 2 &
\centering $a= 0.6557 + 0.0357i$, $b= 0.8491 + 0.9340i$, and $r=-99.554 + 9i$
&
\begin{center}
$(0.6787 + 0.7577i, 0.7431 + 0.3922i, 0.6555 + 0.1712i)$
\end{center}
&
\begin{center}
Convergent orbit and converge to $(0.0109 + 0.0205i,  -0.2251 + 0.0216i, 0.0009 + 0.0150i)$

\end{center}
\\
\hline
\centering 3 &
\centering $a= 0.8147 + 0.9058i$, $b=0.1270 + 0.9134i
$ and $r=-31.0966 - 95i$
&
\begin{center}
$(0.6324 + 0.0975i, 0.2785 + 0.5469i, 0.9575 + 0.9649i)$
\end{center}
&
\begin{center}
Knotted Quasi-periodic Orbit.
\end{center}
\\
\hline
\centering 4 &
\centering $a=0.6557 + 0.0357i$, $b= 0.8491 + 0.9340i$ and $r= 61.3742 -78.9i$
&
\begin{center}
$(0.6787 + 0.7577i, 0.7431 + 0.3922i,  0.6555 + 0.1712i)$
\end{center}
&
\begin{center}
Knotted Quasi-periodic Orbit.
\end{center}
\\
\hline
\centering 5 &
\centering $a=0.0759 + 0.0540i$, $b=0.5308 + 0.7792i$ and $r=65.6265 -75.4i$
&
\begin{center}
$(0.9340 + 0.1299i, 0.5688 + 0.4694i, 0.0119 + 0.3371i)$
\end{center}
&
\begin{center}
Knotted Chaotic orbit.
\end{center}
\\
\hline

\end{tabular}
\caption{Dynamics of discrete complex Lorenz system when $\abs{r}=99.96$ and $dt=0.00005$.}
\label{Table:1}
\end{table}

\noindent
Here for all the five examples sited in the Table $6$, the modulus of $r$ is $99.96$ and it is found that all the $3$ Dimensional trajectories are knotted with different kind of asymptotic behaviour (chaotic, quasi periodic, convergent) as noted in the Table $6$ except the case in serial number $2$. The case in serial number $2$ in the Table $6$ is an example of a trajectory which is convergent to $(0,0,0)$ where $\abs{r}=99.96$.\\

\noindent
It is noted that for the parameters $a, b,$ and $r$ as noted in serial numbers $1$ to $5$ in the Table $6$, it is found that $\abs{r} \geq \abs{a \frac{a+b+3}{a-b-1}}$ which ensures that none of the trajectories are convergent to the non-origin equilibriums.

\begin{figure}[H]
      \centering

      \resizebox{16cm}{!}
      {
      \begin{tabular}{c c c}
      \includegraphics [scale=8]{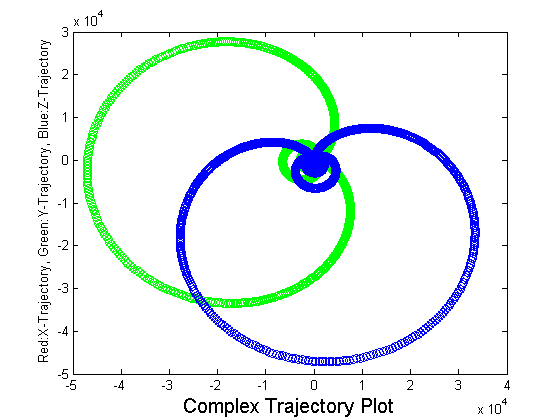}
      \includegraphics [scale=8]{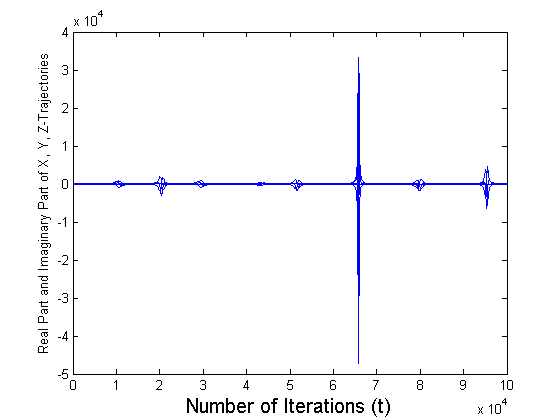}
      \includegraphics [scale=8]{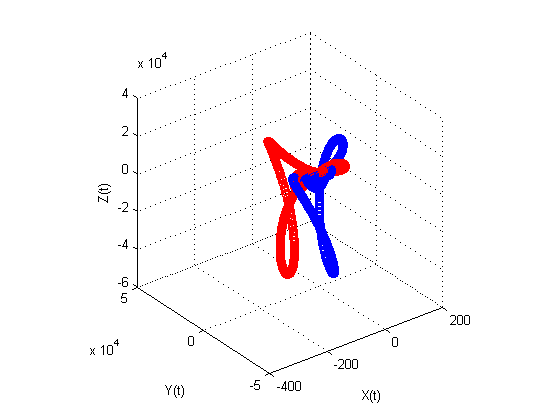}\\
      \includegraphics [scale=8]{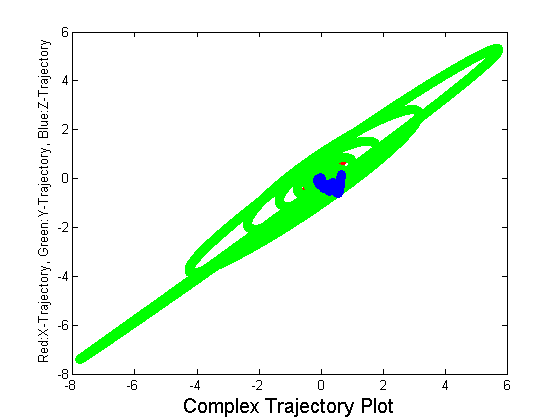}
      \includegraphics [scale=8]{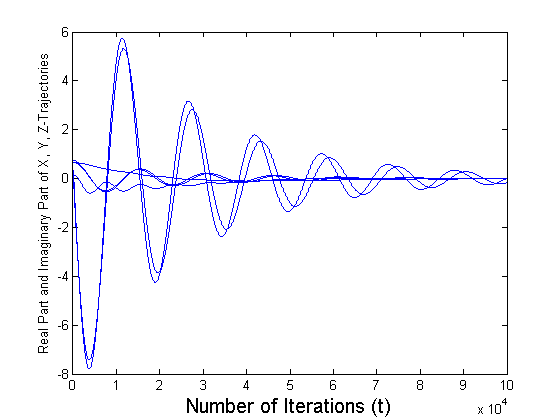}
      \includegraphics [scale=8]{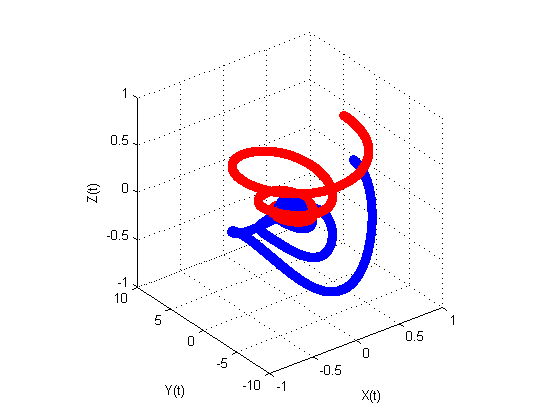}\\
      \includegraphics [scale=8]{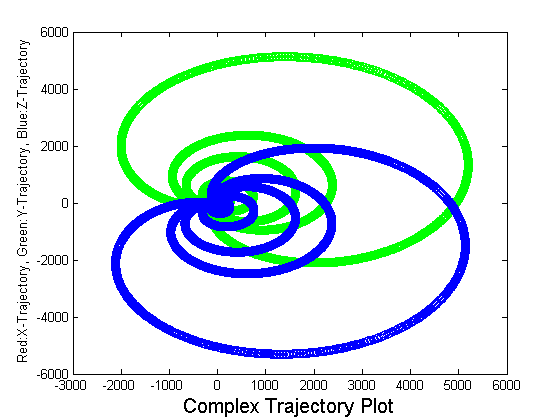}
      \includegraphics [scale=8]{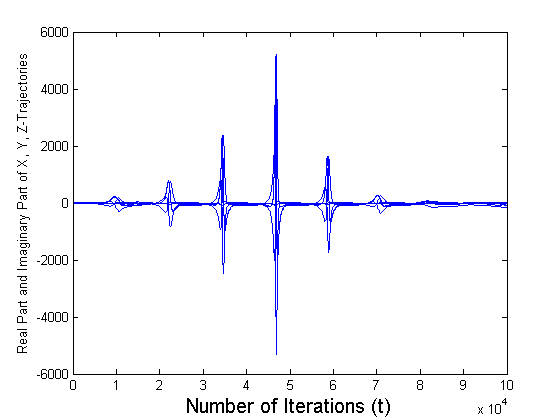}
      \includegraphics [scale=8]{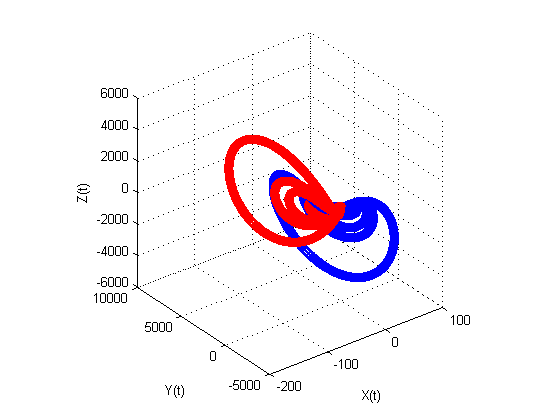}\\
      \includegraphics [scale=8]{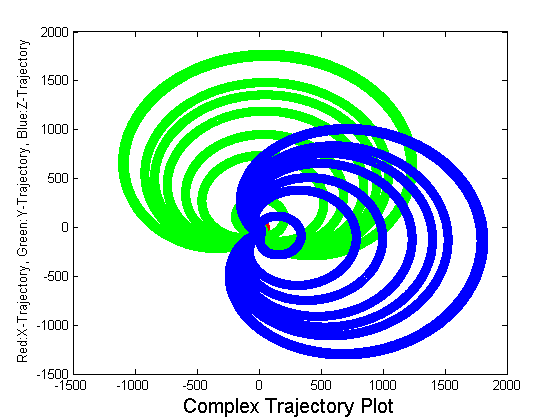}
      \includegraphics [scale=8]{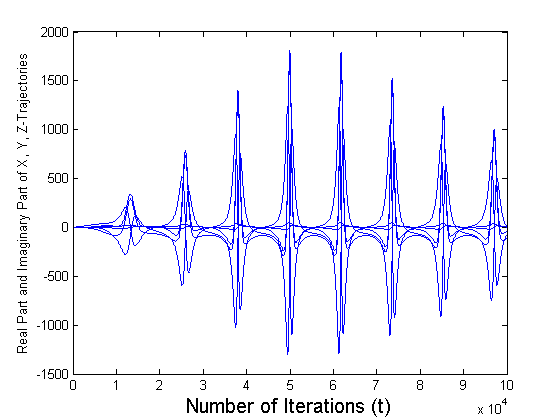}
      \includegraphics [scale=8]{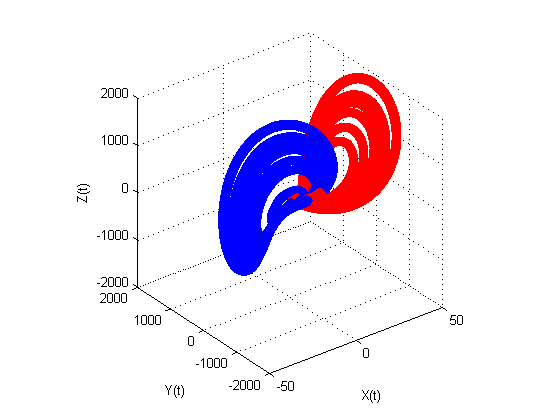}\\
      \includegraphics [scale=8]{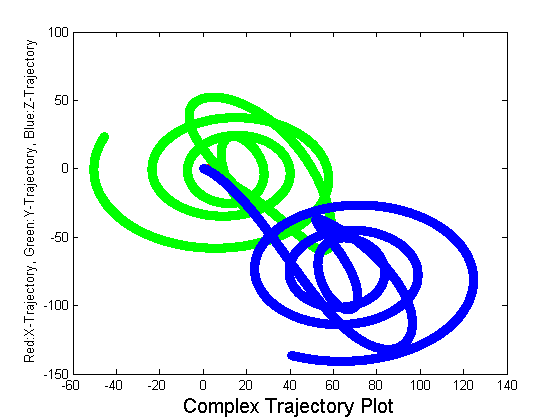}
      \includegraphics [scale=8]{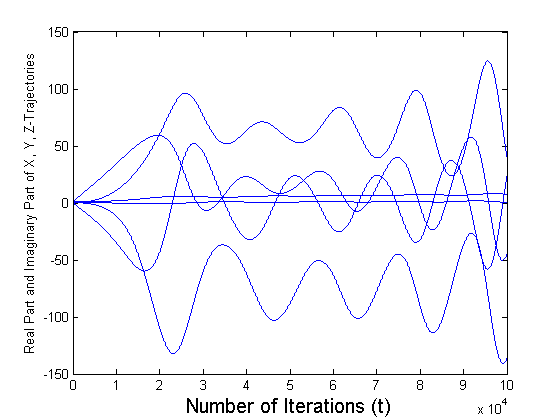}
      \includegraphics [scale=8]{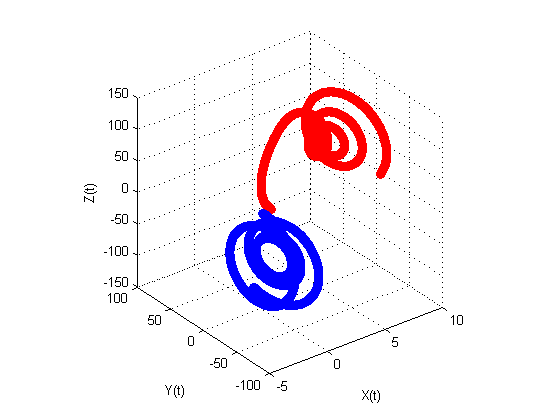}\\

            \end{tabular}
      }
\caption{Serial number $1$ to $5$ in the Table 6, from top to bottom row: Top Left: Complex trajectory plots, Top Middle: X, Y, Z-trajectories, Top Right: Trajectory plot (Red: Real part, Blue: Imaginary part) in 3 dimension.}
      \begin{center}

      \end{center}
      \end{figure}

\noindent
Here in the Table $7$, five examples are taken where five different values of r have been considered to look into the asymptotic behavior of the discrete complex Lorenz system Eq.$(3,4,5)$. All the four examples serial no $1$ to $4$ are chaotic where it is found that largest Lyapunav exponent is positive and the last one in serial number $5$ the trajectory is asymptotically convergent and converging to $(0,0,0)$.

\begin{table}[H]

\begin{tabular}{| m{1cm}||m{3.8cm}| |m{4.5cm}| |m{4.5cm}|}
\hline \centering \textbf{Serial No} &
\begin{center}
\textbf{Parameters: $a$, $b$ and $r$}
\end{center}
 &
\begin{center}
\textbf{Initial Value}
\end{center}
&
\begin{center}
\textbf{Dynamics in the Discrete Complex Lorenz System}
\end{center}
\\
\hline
\centering 1 &
\centering $a= 0.0046 + 0.7749i$, $b=0.8687 + 0.0844i$ and $r=635$
&
 \begin{center}
 $(0.3998 + 0.2599i, 0.8001 + 0.4314i, 0.9106 + 0.1818i)$
 \end{center}
&
\begin{center}
Chaotic (Largest Lyapunav exponent is positive.
\end{center}
\\
\hline
\centering 2 &
\centering $a= 0.9448 + 0.4909i$, $b= 0.3377 + 0.9001i$, and $r=-22$
&
\begin{center}
$(0.3692 + 0.1112i, 0.7803 + 0.3897i, 0.2417 + 0.4039i)$
\end{center}
&
\begin{center}
Chaotic (Largest Lyapunav exponent is positive).

\end{center}
\\
\hline
\centering 3 &
\centering $a= 0.2625 + 0.8010i$, $b=0.9289 + 0.7303i
$ and $r=-195$
&
\begin{center}
$(0.4886 + 0.5785i, 0.2373 + 0.4588i, 0.9631 + 0.5468i)$
\end{center}
&
\begin{center}
Chaotic (Largest Lyapunav exponent is positive)
\end{center}
\\
\hline
\centering 4 &
\centering $a= 0.2625 + 0.8010i$, $b=0.9289 + 0.7303i
$ and $r=195$
&
\begin{center}
$(0.4886 + 0.5785i, 0.2373 + 0.4588i, 0.9631 + 0.5468i)$
\end{center}
&
\begin{center}
Chaotic (Largest Lyapunav exponent is positive)
\end{center}
\\
\hline
\centering 5 &
\centering $a=0.7962 + 0.0987i$, $b=0.3354 + 0.6797i$ and $r=-148$
&
\begin{center}
$(0.1366 + 0.7212i, 0.1068 + 0.6538i, 0.4942 + 0.7791i)$
\end{center}
&
\begin{center}
Asymptotically converges to $(0,0,0)$.
\end{center}
\\
\hline

\end{tabular}
\caption{Dynamics of Discrete Complex Lorenz System when for five different values of $r$ and $dt=0.00005$.}
\label{Table:1}
\end{table}

\begin{figure}[H]
      \centering

      \resizebox{16cm}{!}
      {
      \begin{tabular}{c c c}
      \includegraphics [scale=8]{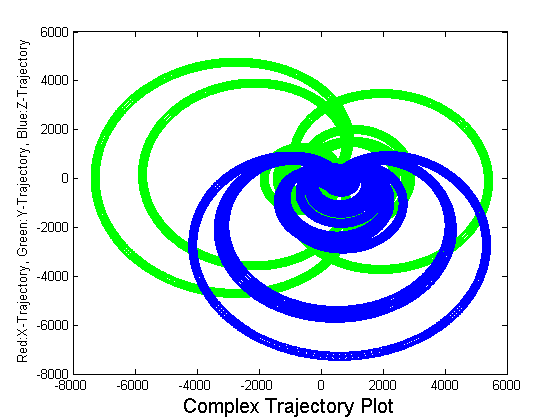}
      \includegraphics [scale=8]{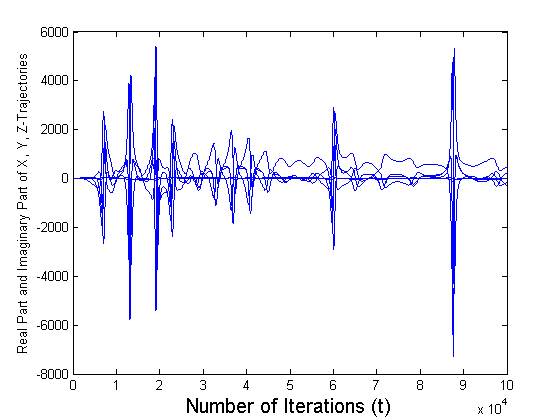}
      \includegraphics [scale=8]{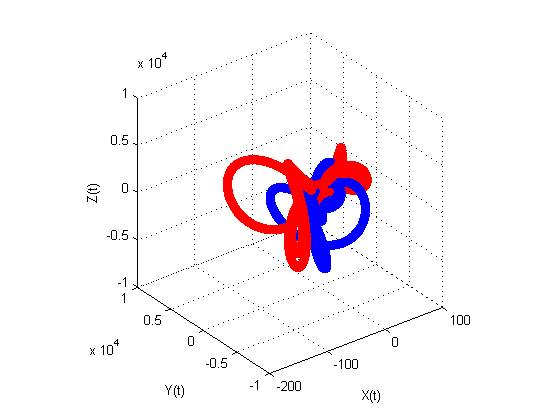}\\
      \includegraphics [scale=8]{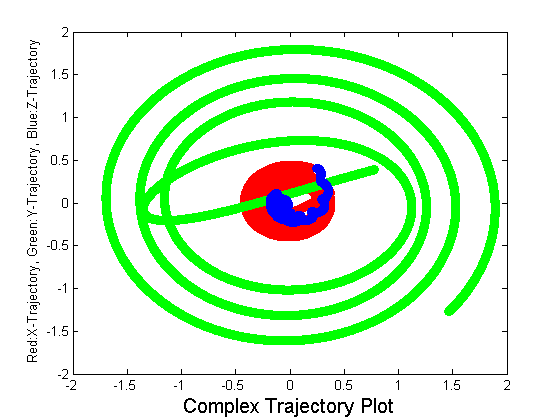}
      \includegraphics [scale=8]{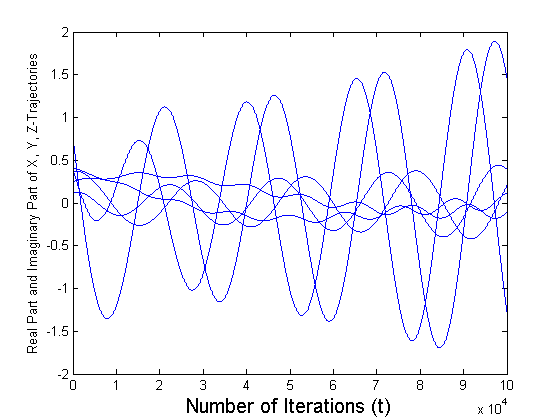}
      \includegraphics [scale=8]{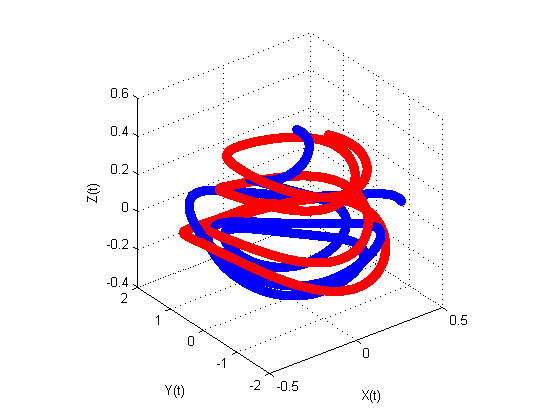}\\
      \includegraphics [scale=8]{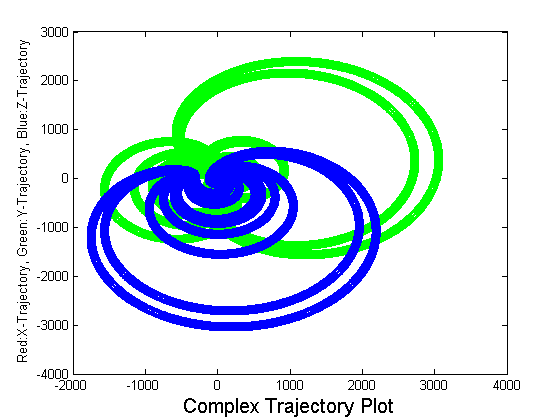}
      \includegraphics [scale=8]{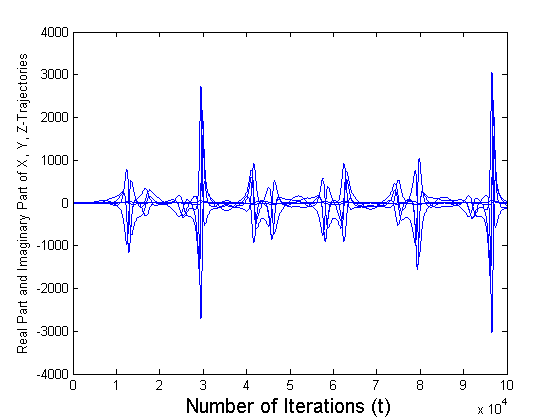}
      \includegraphics [scale=8]{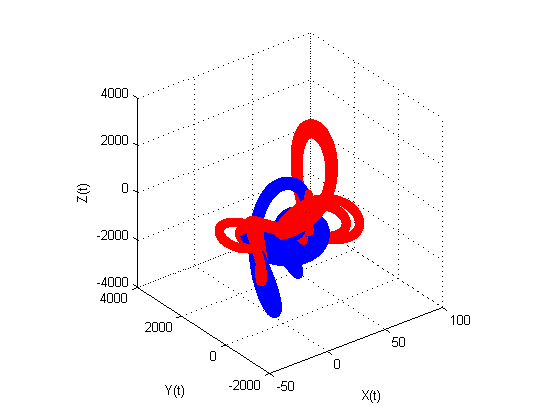}\\
      \includegraphics [scale=8]{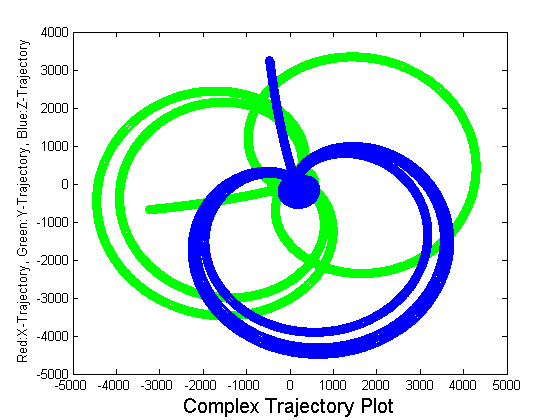}
      \includegraphics [scale=8]{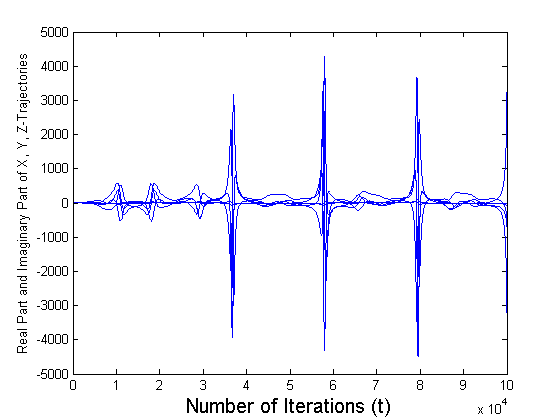}
      \includegraphics [scale=8]{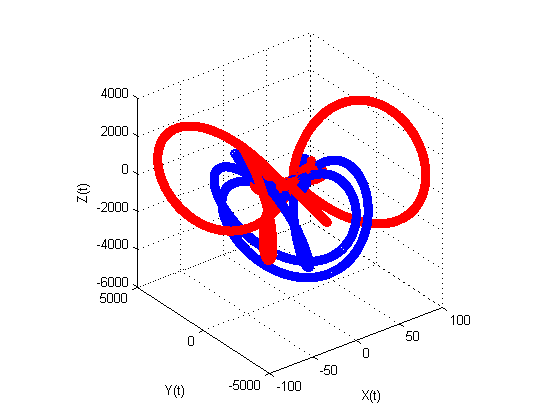}\\
      \includegraphics [scale=8]{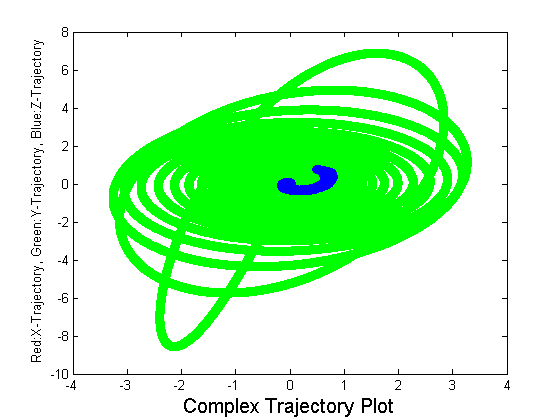}
      \includegraphics [scale=8]{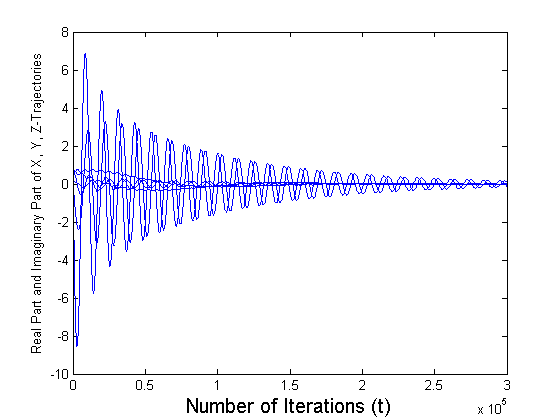}
      \includegraphics [scale=8]{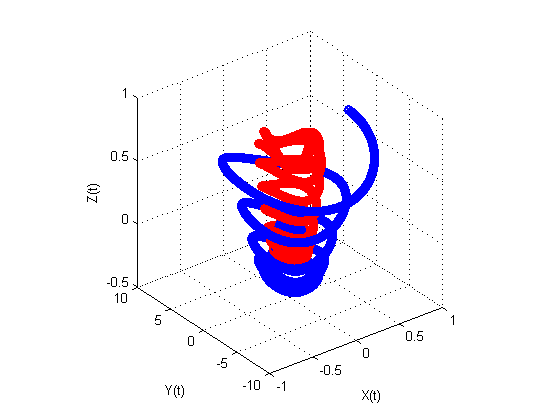}\\

            \end{tabular}
      }
\caption{Serial number $1$ to $5$ in the Table 7, from top to bottom row: Top Left: Complex trajectory plots, Top Middle: X, Y, Z-trajectories, Top Right: Trajectory plot (Red: Real part, Blue: Imaginary part) in 3 dimension.}
      \begin{center}

      \end{center}
      \end{figure}

\noindent
Here we shall encounter a few examples where one of the parameters either $a$ or $b$ is zero. It is noted that when $a=0$, then from Eq. $(3,4,5)$, it comes out that $x_{k+1}=x_k$ and other two equations Eq.$(4,5)$ remain same. When $b=0$, the Eq. (5) becomes $z_{k+1}=z_k+x_ky_kdt$ and other two equations Eq.$(3,4)$ remain same.

\begin{table}[H]

\begin{tabular}{| m{1cm}||m{3.8cm}| |m{4.5cm}| |m{4.5cm}|}
\hline \centering \textbf{Serial No} &
\begin{center}
\textbf{Parameters: $a$, $b$ and $r$}
\end{center}
 &
\begin{center}
\textbf{Initial Value}
\end{center}
&
\begin{center}
\textbf{Dynamics in the Discrete Complex Lorenz System}
\end{center}
\\
\hline
\centering 1 &
\centering $b=0.8491 + 0.9340i$ and $r=7+i$
&
 \begin{center}
 $(0.6787 + 0.7577i, 0.7431 + 0.3922i, 0.6555 + 0.1712i)$
 \end{center}
&
\begin{center}
Converges to  $(0.6787 + 0.7577i, 3.4194 + 2.3119i, 3.0034 + 1.7374i)$
\end{center}
\\
\hline
\centering 2 &
\centering $b= 0.1361 + 0.8693i$, and $r=-218-255i$
&
\begin{center}
$(0.5797 + 0.5499i, 0.1450 + 0.8530i, 0.6221 + 0.3510i)$
\end{center}
&
\begin{center}
Chaotic (Largest Lyapunav exponent is positive).

\end{center}
\\
\hline
\centering 3 &
\centering $a=0.9448 + 0.4909i$, and $r=-1-17i$
&
\begin{center}
$(0.9001 + 0.3692i, 0.1112 + 0.7803i, 0.3897 + 0.2417i)$
\end{center}
&
\begin{center}
Converges to $(0.0036i, -0.0143 + 0.0081i, 3.1585 -23.6625i)$
\end{center}
\\
\hline
\centering 4 &
\centering  $a=0.4039 + 0.0965i$ and $r=-25+8i$
&
\begin{center}
$(0.9561 + 0.5752i, 0.0598 + 0.2348i, 0.3532 + 0.8212i)$
\end{center}
&
\begin{center}
Converges to $(0.0020 - 0.0012i, 0.0142 + 0.0106i, -0.3889 + 0.1130i)$
\end{center}
\\
\hline
\centering 5 &
\centering $a=0.0154 + 0.0430i$ and $r=-24-4i$
&
\begin{center}
$(0.7317 + 0.6477i, 0.4509 + 0.5470i, 0.2963 + 0.7447i)$
\end{center}
&
\begin{center}
Transient Chaos
\end{center}
\\
\hline

\end{tabular}
\caption{Dynamics of Discrete Complex Lorenz System when for five different values of $r, b$ where either $a=0$ or $b=0$ and $dt=0.00005$.}
\label{Table:8}
\end{table}

\noindent
The fixed points for the parameters (Refer Serial no. 1) of the system Eq.(3,4,5) are $(\pm 2.83157 + 1.46935i, \pm 2.83157 + 1.46935i, 8 + i)$. It is noted that $\abs{r}=7.07107$ and $\abs{a\frac{a+b+3}{a-b-1}}=0$, hence the condition for convergence $\abs{r}<\abs{a\frac{a+b+3}{a-b-1}}$ to the fixed points is violating. It is seen that the trajectory is convergent and converges to $(0.6787 + 0.7577i, 3.4194 + 2.3119i, 3.0034 + 1.7374i)$. Here the convergence depends on the initial condition unlike other cases.\\

\noindent
The fixed points for the parameters (Refer Serial no. 3) of the system Eq.(3,4,5) is $(0,0,-17i)$. It is noted that $\abs{r}=17.0294$ and $\abs{a\frac{a+b+3}{a-b-1}}=8.5679$, hence the condition for convergence $\abs{r}<\abs{a\frac{a+b+3}{a-b-1}}$ to the fixed point is not valid. it is found that the trajectory in this case is convergent and converges to $(0.0036i, -0.0143 + 0.0081i, 3.1585 -23.6625i)$.

\noindent
The fixed points for the parameters (Refer Serial no. 4) of the system Eq.(3,4,5) is $(0,0,-24+8i)$. It is noted that $\abs{r}=26.2488$ and $\abs{a\frac{a+b+3}{a-b-1}}=2.3418$, hence the condition for convergence $\abs{r}<\abs{a\frac{a+b+3}{a-b-1}}$ to the fixed point is invalid. it is found that the trajectory in this case is convergent and converges to $(0.0020 - 0.0012i, 0.0142 + 0.0106i, -0.3889 + 0.1130i)$.

\noindent
For the parameters (Refer Serial no. 5), the trajectory possesses to transient chaos where one positive Lyapunav exponent is found.

\begin{figure}[H]
      \centering

      \resizebox{16cm}{!}
      {
      \begin{tabular}{c c c c}
      \includegraphics [scale=8]{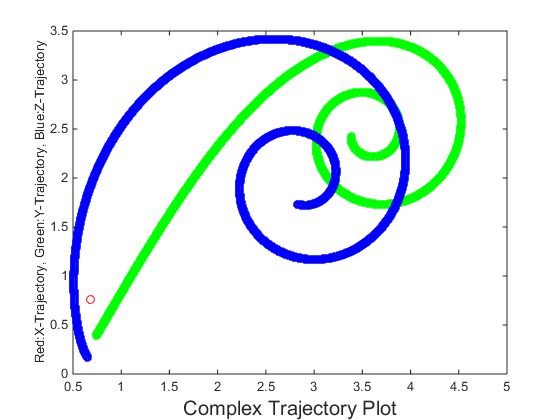}
      \includegraphics [scale=8]{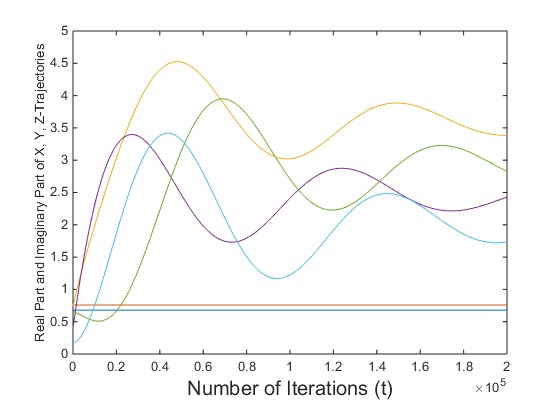}
      \includegraphics [scale=8]{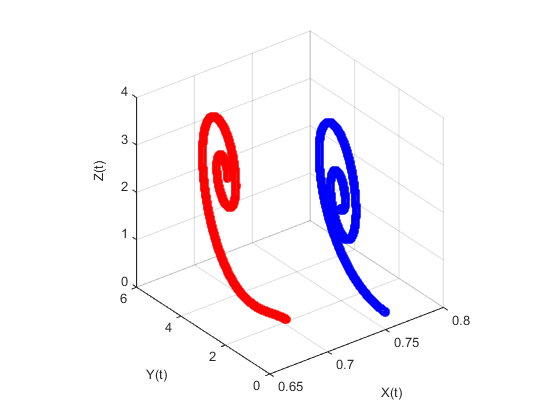}\\
      \includegraphics [scale=8]{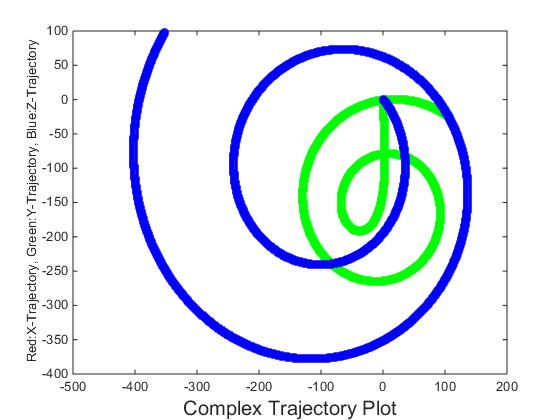}
      \includegraphics [scale=8]{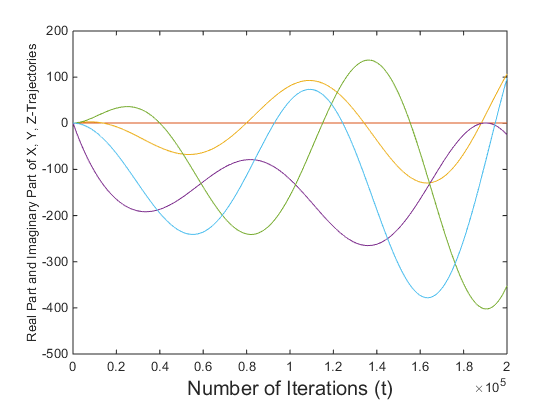}
      \includegraphics [scale=8]{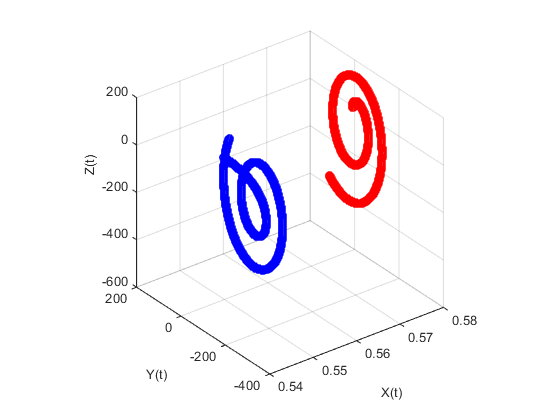}\\
      \includegraphics [scale=8]{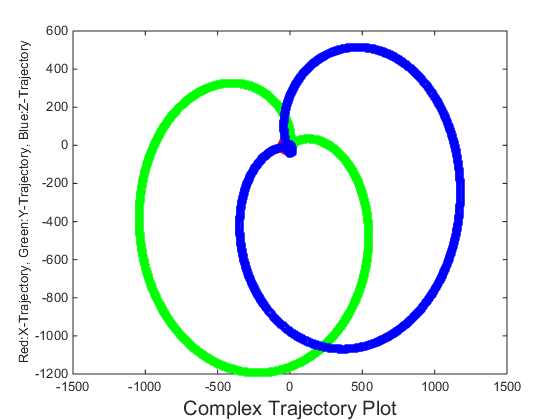}
      \includegraphics [scale=8]{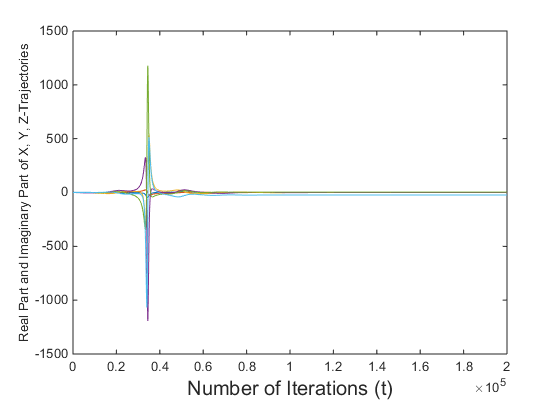}
      \includegraphics [scale=8]{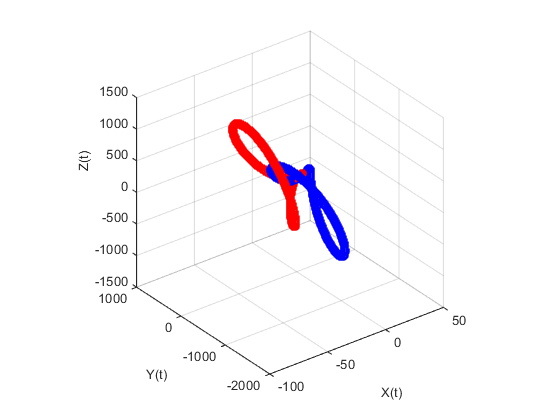}\\
      \includegraphics [scale=8]{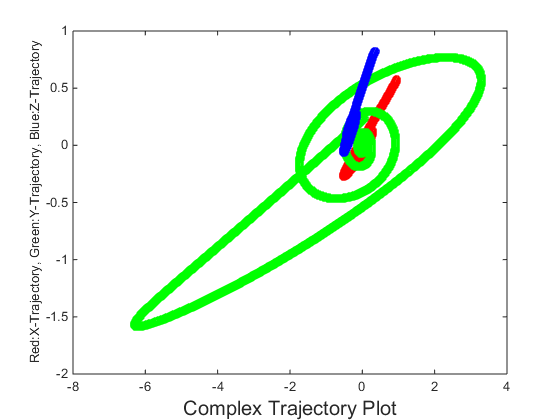}
      \includegraphics [scale=8]{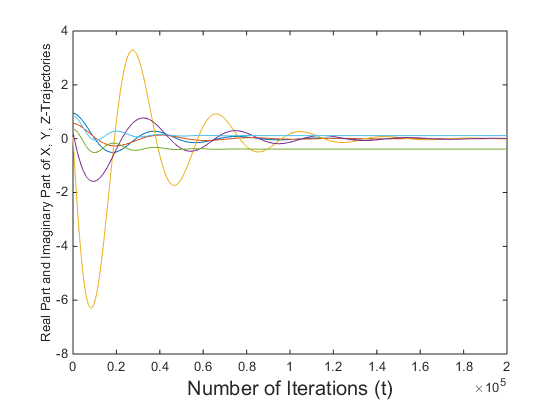}
      \includegraphics [scale=8]{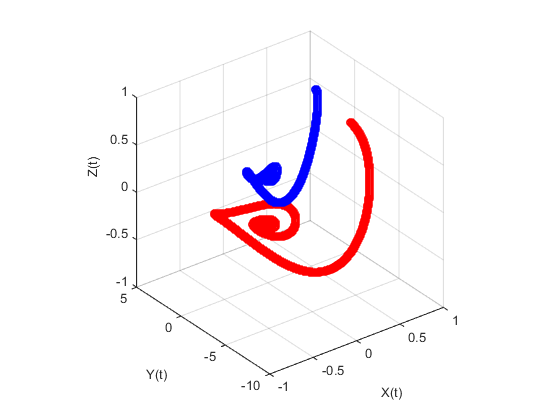}\\
      \includegraphics [scale=8]{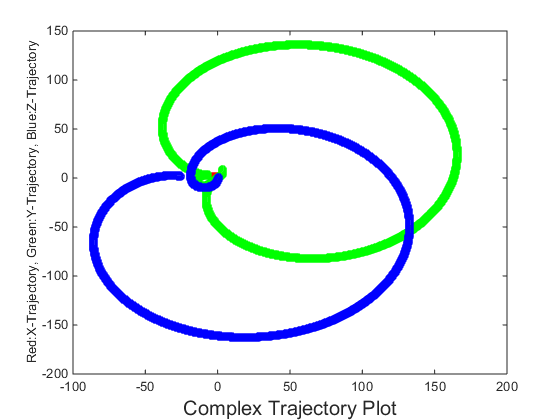}
      \includegraphics [scale=8]{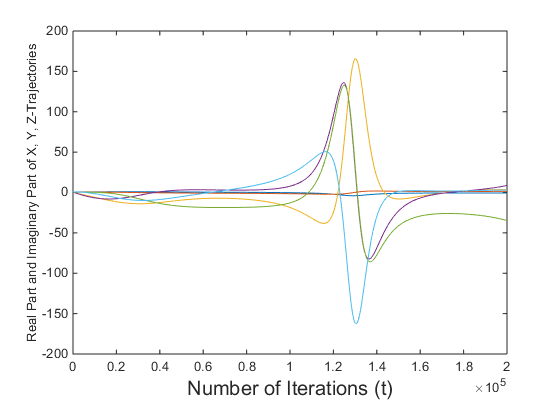}
      \includegraphics [scale=8]{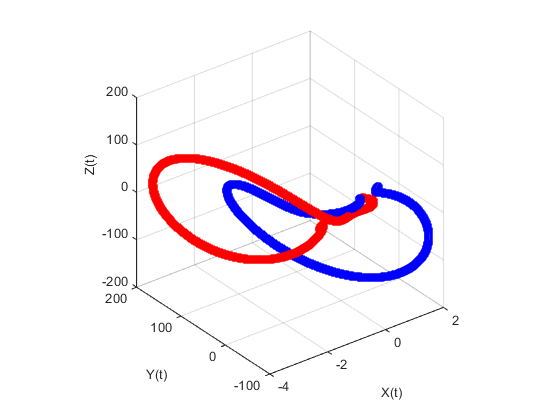}\\

            \end{tabular}
      }
\caption{Serial number $1$ to $5$ in the Table 8, from top to bottom row: Top Left: Complex trajectory plots, Top Middle: X, Y, Z-trajectories, Top Right: Trajectory plot (Red: Real part, Blue: Imaginary part) in 3 dimension.}
      \begin{center}

      \end{center}
      \end{figure}

\section*{Acknowledgement}
The author thanks to \emph{Snighdha Das of IIT Kharagpur and Prof. Parichaoy Kumar Das of UPES, Dehradun} for long hours fruitful discussions and suggestions.


\end{document}